\documentclass[secthm,seceqn,amsthm,ussrhead]{amsart}
\usepackage{amsmath,latexsym}
\usepackage[psamsfonts]{amssymb}
\usepackage{times}
\usepackage[mathcal]{euscript}
\numberwithin{equation}{section}
\usepackage{enumerate}

\newcommand{\bbT}{\mathbb T}

\renewcommand{\epsilon}{\varepsilon}
\newcommand{\be}{\begin{equation}}
\newcommand{\ee}{\end{equation}}
\newcommand{\no}{\nonumber}
\newcommand{\F}{\mathbb{F}}

\newcommand{\N}{\mathbb{N}}

\newcommand{\R}{\mathbb{R}}
\renewcommand{\S}{\mathbb{S}}
\newcommand{\T}{\mathbb{T}}

\newcommand{\Z}{\mathbb{Z}}
\newcommand{\cE}{{\mathcal E}}
\newcommand{\cF}{{\mathcal F}}

\newcommand{\cH}{{\mathcal H}}

\newcommand{\cU}{{\mathcal U}}

\DeclareMathOperator{\Ker}{\mathrm{Ker}}

{\bf}{\it}
\newtheorem{theorem}{Theorem}[section]
\newtheorem{lemma}[theorem]{Lemma}

\newtheorem{hypothesis}[theorem]{Hypothesis}
\newtheorem{definition}[theorem]{Definition}

\newtheorem{remark}[theorem]{Remark}

\date{\today}

\begin{document}
\title[Discrete spectrum asymptotics for the three-particle
Hamiltonians ...]
 {Discrete spectrum asymptotics for the three-particle
Hamiltonians on lattices}

%\title{ Asymptotics for the number of eigenvalues of a three-particle
%Schr\"{o}dinger operator on lattice:(the case of two identical
%massive fermions)\\}

 \author{Sergio Albeverio$^{1,2,3}$,Saidakhmat
 N. Lakaev $^{4,5}$ Axmad M. Xalxo'jaev $^5$}

%\author{S.  Albeverio$^{1,2,3}$}
\address{$^1$ Institut f\"{u}r Angewandte Mathematik,
Universit\"{a}t Bonn (Germany)} \email{albeverio@uni.bonn.de}
\address{
$^2$ \ SFB 611, \ Bonn, \ BiBoS, Bielefeld - Bonn\ (Germany)}
\address{
$^3$  CERFIM, Locarno and USI (Switzerland)}

%\author{S.  N. Lakaev$^{5,6}$}
\address{
$^4$ Samarkand Division of Academy of sciences of Uzbekistan
(Uzbekistan)} \email{lakaev@yahoo.com}

\address{
$^5$ Samarkand State University,Samarkand (Uzbekistan)} \email
{ahmadx@mail.ru}

\begin{abstract}
We consider the Hamiltonian of a system of three quantum mechanical
particles on the three-dimensional lattice $\Z^3$ interacting via
short-range pair potentials.

We prove for the two-particle energy operator $h(k),$ $k\in \T^3$
the two-particle quasi-momentum, the existence of a unique positive
eigenvalue $z(k)$ lying below the essential spectrum under
assumption that the operator $h(0)$ corresponding to the zero value
of $k$ has a zero energy resonance.

We describe the location of the essential spectrum of the
three-particle discrete Schr\"{o}dinger operators $H(K)$,$K$ the
three-particle quasi-momentum by the spectra of $h(k),\,k\in \T^3.$

We prove the existence of infinitely many eigenvalues of $H(0)$ and
establish  for the number of eigenvalues $N(0,z)$ lying below $z<0$
the asymptotics

\begin{equation*}\label{asimz} \lim\limits_{z \to
-0}\frac{N(0,z)}{|\log |z||}=\frac{\lambda _0}{2\pi},
\end{equation*}
where   $\lambda _0$ a unique positive solution of the equation
$$ \lambda = \frac{8 \sinh \pi\lambda /6}{\sqrt 3 \cosh \pi\lambda/2}.$$

We prove that for all $ K \in U_\delta^0(0),$ where $U_\delta^0(0)$
some punctured $\delta
>0$ neighborhood of the origin, the number $N(K,0)$ of
eigenvalues  the operator $H(K)$ below zero is finite and satisfy
the asymptotics
\begin{equation*}\label{asimk}
\lim\limits_{|K| \to 0}\frac{N(K,0)}{|\log |K||}=\frac{\lambda
_0}{\pi}.
\end{equation*}

\end{abstract}
\maketitle



\section{ Introduction}
We consider a system of three identical particles (bosons) on the
three-dimensional lattice $\Z^3$ interacting by means of short-range
pair attractive potentials.

The main goal of the present paper is to prove the finiteness or
infiniteness (Efimov's effect) of the number of eigenvalues lying
below zero (the bottom of the essential spectrum of $H(0)$) of the
three-particle discrete Schr\"{o}dinger operator $H(K)$ depending on
the total quasi-momentum $ K \in U_\delta^0(0),$ where $U_\delta(0)$
a $\delta
>0$ neighborhood of the origin.

Efimov`s  effect is one of the remarkable results in the spectral
analysis for   three-particle Schr\"{o}dinger operators associated
to a system of three particles moving on Euclid space $\R^3$ or
integer lattice $\Z^3$: if none of the three two-particle
Schr\"{o}dinger operators (corresponding to the two-particle
subsystems) has negative eigenvalues, but at least two of them have
a zero energy resonance, then this three-particle Schr\"{o}dinger
operator has an infinite number of discrete eigenvalues,
accumulating at zero.

Since its discovery by Efimov in 1970 \cite{Efi} much research have
been devoted to this subject. See, for example
\cite{AHW,AmNo,DFT,FaMe,OvSi,Sob,Tam91,Tam94,Yaf74}.

The main result obtained by  Sobolev \cite{Sob} (see also
\cite{Tam94}) is an asymptotics of the form $\cU_0|log|\lambda||$
for the number of eigenvalues below $\lambda,\lambda<0$, where the
coefficient ${\cU}_0$ does not depend on the two-particle potentials
$ v_\alpha $ and is a positive function of the ratios
$m_1/m_2,m_2/m_3$ of the masses of the three-particles.

Recently  the existence of Efimov's effect for $N$-body quantum
systems with $N\geq 4$ has been proved by X.P. Wang in \cite{Wang}.

In fact in \cite{Wang} for the total (reduced) Hamiltonian a lower
bound on the number of eigenvalues of the form
$C_0|log(E_0-\lambda)|$ is given, when $\lambda$ tends to $E_0$,
where $C_0$ is a positive constant and $E_0$ is the bottom of the
essential spectrum.

 The kinematics of quantum particles
on lattices, even in the two and three-particle sector, is rather
exotic. For instance,
 due to the fact that the discrete analogue of
the Laplacian or its generalizations  are not
translationally invariant, the Hamiltonian of a system  does not
separate into two parts, one relating to the center-of-mass motion
and the other one relating to the internal degrees of freedom.

As a consequence any local substitute of the effective mass-tensor
(of a ground state) depends on the quasi-momentum of the system and,
in addition, it  is only semi-additive (with respect to the partial
order on the set of positive definite matrices). This is the
so-called {\it excess mass} phenomenon for lattice systems (see,
e.g., \cite{Mat} and  \cite{Mog91}): the effective mass of the bound
state of an $N$-particle system is greater than and in general,
 not equal to the sum of the effective
masses of the constituent quasi-particles.

The three-body  problem on lattices  can  be reduced to the
effective three-particle Schr\"odinger operators  by using the
Gelfand transform. The underlying Hilbert space $\ell^2((\Z^d)^3)$
 is decomposed as a direct von Neumann integral associated with the
representation of the discrete group $\Z^3$  by shift operators on
the lattice and  the total three-body Hamiltonian appears to be
decomposable. In contrast to the continuous case, the corresponding
fiber Hamiltonians $H(K) $ associated with  the direct decomposition
depend parametrically  the quasi-momentum,  $K\in
\T^3=(-\pi,\pi]^3$, which  ranges over a cell of the dual lattice.
Due to the loss of the spherical symmetry of the problem,  the
spectra of the family $H(K)$ turn out to be  rather sensitive to the
quasi-momentum $K \in \T^3.$

In particular,  Efimov's effect exists only for the zero value of
the three-particle quasi-momentum $K,$ which is proven for
Hamiltonians of a system of three particles interacting via pair
zero-range attractive potentials on $\Z^3$(see, e.g.,
\cite{ALM98,ALzM04,Ltmf91,Lfa93, LaAb99,Mat} for relevant
discussions and \cite{ALMM,FIC,GrSc,KM,Mat,MS,Mog91,RSIII} for the
general study of the low-lying excitation spectrum for quantum
systems on lattices).

 Denote by $\tau (K)$ the bottom of the essential
spectrum of the three-particle discrete Schr\"{o}dinger operator
$H(K),\,K\in \T^{3}$ and by $N(K,z)$ the number of eigenvalues lying
below $z\leq \tau (K).$

The main results of the present paper are as follows:

(i) the operator $H(0)$ has infinitely many eigenvalues below the
bottom of the essential spectrum and for the number of eigenvalues
$N(0,z )$ lying below $z<0$ the asymptotics
\begin{equation*}\label{asimz} \lim\limits_{z \to
-0}\frac{N(0,z)}{|\log |z||}=\frac{\lambda _0}{2\pi},
\end{equation*}
holds, where   $\lambda _0$ a unique positive solution of the equation
$$ \lambda = \frac{8 \sinh \pi\lambda /6}{\sqrt 3 \cosh \pi\lambda/2}.$$
This result is similar to the asimptotics founded in the continuous
case by Sobolev \cite{Sob};

(ii) for some  punctured $\delta
>0$ neighborhood $U_\delta^0(0)$ of the origin and for all $ K \in U_\delta^0(0)$ the
 number $N(K,0)$ is a finite and satisfy the following
asymptotics
\begin{equation*}\label{asimk} \lim\limits_{|K| \to
0}\frac{N(K,0)}{|\log |K||}=\frac{\lambda _0}{\pi}.
\end{equation*}
This result is characteristic for the lattice system and does not
have any analogue in the continuous case .

We underline   that these results are  in contrast to similar results for the
continuous three particle Schr\"odinger operators, where the number
of eigenvalues does not depend on the three-particle total momentum
$K \in R^3.$

Moreover these results  are also in contrast with the results for
two-particle operators, in which  discrete Schr\"odinger operators have
finitely many eigenvalues for
 all $k \in U_\delta(0),$ where $U_\delta(0)=\{K \in \T^3: |K|<\delta\}$
 is  a $\delta-$ neighborhood of the origin.

Note that to prove these results in the present paper we derive  an
asymptotics of the Birman-Schwinger operator $G(k,0)$ resp. $G(0,z)$
as $k\rightarrow0$ resp. $z\rightarrow0.$ In particular, we prove
that the operator valued function $G(k,0)$ resp. $G(0,z)$ is
differentiable in $|k|$ at $k=0 \in T^3$ resp.in $z$ at $z=0.$

This result  has been proved in the continuum case (see
\cite{Sob,Tam94}) using resolvent expansion established in
\cite{JeKa}.

The organization of the present paper is as follows.

Section 1 is an introduction.

In Section 2 we introduce the Hamiltonians of systems of two and
three-particles in coordinate and momentum representations as
bounded self-adjoint operators in the corresponding Hilbert spaces.

In Section 3 we introduce the total quasi-momentum and decompose the
energy operators into von Neumann direct integrals, choosing
relative coordinate systems.

In section 4 we introduce the concept of a zero energy resonance
(threshold resonance).
 In Section 5 we state the
main results of the paper.

In Section 6 we study spectral properties of the two-particle
discrete Schr\"{o}dinger operator $h(k),k\in \T^3.$ We prove the
existence of a positive eigenvalue below the bottom of the essential
spectrum of $h(k),k\in \T^3$ (Theorem \ref{mavjud}) and obtain an
asymptotics for the Birman-Schwinger operators associated to
$h(k),k\in \T^3.$

In Section 7 we introduce the {\it channel operator} and describe
its spectrum by the spectrum of the two-particle discrete
Schr\"{o}dinger operators (Theorem \ref{ess}).

In Section 8 we prove  the Birman-Schwinger principle for the three
identical particle Schr\"{o}dinger operator on lattice $\Z^3$.

In Section 9 we  derive the asymptotics for the number of eigenvalues
$N(0,z)$ resp. $N(K,0)$ of $H(0)$  as $z\rightarrow 0$ resp. $H(K)$
as $|K|\rightarrow 0$ (Theorem \ref{asimpZ} resp.Theorem
\ref{asimpK}).

Throughout the present paper we adopt the following notations: We
denote by  $\T^3$   the three-dimensional torus, i.e.,the cube
$(-\pi,\pi]^3$ with appropriately  identified sides. The torus
$\T^3$ will always be considered as an abelian group with respect to
the addition and multiplication by real numbers regarded as
operations on the three-dimensional space $\R^3$ modulo $(2\pi
\Z)^3$.

 For each (sufficiently small) $\delta>0$ the notation
 $ U_\delta(0)=\{K \in \T^3: |K|<\delta\}$
 stands for a $\delta-$ neighborhood of the origin and
 $U^0_\delta(0)=U_\delta(0)\setminus \{0\} $ for a punctured
$\delta-$  neighborhood. The subscript $ \alpha $ (and also $\beta $
and $\gamma $) always runs from 1 to 3 and we use the convention  $
\alpha \not=\beta, \beta\not= \gamma , \gamma \not=\alpha  $.

\section{ Description of the energy operators of two and three arbitrary
particles on a lattice and formulations of the main results}

Let ${\Z}^{3}$ be the three-dimensional lattice and let $
\,\,(\Z^{3})^{m}\,,\,\,m\in N$ be the Cartesian $m-$ th  power of $
\,\,\Z^{3}.\,\,$ Denote by $\ell _{2}((\Z^{3})^{m})$ the Hilbert
space of square-summable functions $\,\,{\varphi}$ defined on
$(\Z^{3})^{m}$ and let  $\ell _{2}^{s}((\Z^{3})^{m})\subset \ell
_{2}((\Z^{3})^{m})\,$ be   the subspace of  symmetric functions.

The free Hamiltonian $ \hat h_0$ of a system of
 two identical quantum mechanical particles on the  three dimensional
 lattice $\Z^3$ is defined  by
\begin{equation*}\label{free}
 (\hat {h}^0{\hat\psi})(x_\beta,x_\gamma) =\frac{1}{2}\sum_{|s|=1}
[2{\hat \psi}(x_\beta,x_\gamma)-
{\hat\psi}(x_\beta+s,x_\gamma)-{\hat\psi}(x_\beta,x_\gamma+s)],\quad
{\hat\psi} \in \ell^{(s)}_2(({\Z}^3)^2).
\end{equation*}

 The  Hamiltonian $ \hat h$  of a system of two quantum-mechanical
 identical particles
 interacting through a short-range pair potential $\hat v$
 is usually associated with the following
bounded self-adjoint operator on the Hilbert space $\ell^{(s)}_2((
{\Z}^3)^2)$ and has form
\begin{equation}\label{two-part}
\hat h =\hat h^0-\hat v,
\end{equation}
where
\begin{equation*}
(\hat v \hat \psi)(x_\beta,x_\gamma) =\hat v(x_\beta-x_\gamma)
{\hat\psi}(x_\beta,x_\gamma),\quad {\hat\psi} \in
\ell^{(s)}_2(({\Z}^3)^2).
\end{equation*}

  The free Hamiltonian
$ \widehat H_0$ of a system of
 three identical quantum mechanical particles on the  three-dimensional
 lattice $\Z^3$ is defined  by
\begin{align*}&(\widehat{H}_0\hat{\psi})(x_1,x_2,x_3)\\&=
 \frac{1}{2}\sum_{|s|=1}[3\hat{\psi}(x_1,x_2,x_3)
-\hat{\psi}(x_1+s,x_2,x_3)-\hat{\psi}(x_1,x_2+s,x_3)-
\hat{\psi}(x_1,x_2,x_3+s)] ,
\end{align*}
$$\hat\psi\in \ell^{(s)}_2(({\Z}^3)^3).$$

The  Hamiltonian $ \widehat H$  of a system of three
quantum-mechanical identical particles
 with the two-particle  interaction
$\hat v=\hat v_{\alpha}=\hat
v_{\beta\gamma},\alpha,\beta,\gamma=1,2,3$ is a bounded perturbation
of the free Hamiltonian $ H_0$
\begin{equation}\label{total}
 \widehat {H}=\widehat {H}_0-\widehat {V}_{1}-\widehat {V}_{2}-
\widehat {V}_{3},
\end{equation}
where $\widehat {V}_{\alpha}=\widehat {V}, \alpha=1,2,3 $ is
multiplication operator on $\ell^{(s)}_2(({\Z}^3)^3)$ defined by
\begin{equation*}
(\widehat {V}\hat\psi)(x_1,x_2,x_3)=\hat
{v}(x_\beta -x_\gamma )\hat\psi(x_1,x_2,x_3),
\hat\psi\in \ell^{(s)}_2(({\Z}^3)^3).
\end{equation*}

\begin{hypothesis}\label{hypo}
The function  $\hat v(s)$ is a real, even, nonnegative function on
$\Z^3$ and verify
\begin{equation*}\lim _{|s|\to
 \infty}|s|^{3+\theta} \hat v(s)=0,\quad \theta>\frac{1}{2}.
 \end{equation*}
\end{hypothesis}

\begin{remark}
It is clear that under  of Hypothesis
\ref{hypo}
 the  two-  resp. three- particle Hamiltonian
  \eqref{two-part} resp.  \eqref{total} is a bounded self-adjoint operator on the Hilbert
space  $\ell^{(s)}_2(({\Z}^3)^2)$ resp. $\ell^{(s)}_2(({\Z}^3)^3).$
 \end{remark}
\begin{remark}\label{techn}
We note that  Hypothesis
\ref{hypo}  are far from being optimal,
but we will not discuss further  this point here.
\end{remark}

\subsection{ The momentum  representation}
 Let  ${({\T}^3)^m},\,m\in \bf {N}$ be the  Cartesian
$m$-th power of the torus  ${\T}^3=(-\pi,\pi]^3$ and let
$L^{(s)}_2((\mathbb{T}^{3})^{m})\subset L_2((\mathbb{T}^{3})^{m})$ be
the subspace of symmetric functions.

 Let ${ \mathcal{F}}_m:L_2(({\T}^3)^m) \rightarrow \ell_2((
{\Z}^3)^m)$ be  the standard Fourier transform. Since the subspace
$L^{(s)}_2((\T^{3})^{m})$ is invariant with respect to the group ${
\mathcal{F}}_m,$ i.e.,  $\cF_m L^{(s)}_2((\T^{3})^{m}) \subset
\ell^{(s)} _{2}((\Z^{3})^{m}),$ we denote by  $\cF_m^{s}$ the
restriction of $\cF_m$ on  to  the subspace $L^{(s)}_2((\T^{3})^{m})$.
We   easily  check that $$\cF_m^{s}: L^{(s)}_2((\T^{3})^{m})
\rightarrow \ell^{(s)} _{2}((\Z^{3})^{m}).$$

The two-resp. three-particle Hamiltonians in the momentum
representation are given by bounded self-adjoint operators on the
Hilbert spaces
 $L^{(s)}_2(({\T}^3)^2)$ resp. $L^{(s)}_2(({\T}^3)^3)$
 as follows
 \begin{equation*}
 h={(\cF_2^{s})}^{-1} \hat h \cF_2^{s},\quad
\end{equation*}
 resp.
\begin{equation*}
   H={(\cF_3^{s})}^{-1} \widehat H {\cF_3^{s}}.
 \end{equation*}
The two-particle Hamiltonian $ h$ is of  the form
\begin{equation*}
h = h ^0- v.
 \end{equation*}
 The operator  ${h}^0$ is the multiplication operator by the function
$ \varepsilon (k_1)+\varepsilon (k_2):$
\begin{equation*}
( h^0f)(k_1 ,k_2)=(\varepsilon(k_1)+\varepsilon (k_2 ))f(k_1 ,k_2),\,\,f\in L^{(s)}_2(({\T}^3)^2),
\end{equation*}
where $k_\alpha,\alpha=1,2$ is the {\it quasi-momentum} of the particle
$\alpha.$

The  integral operator  $ v$ is  of convolution type
\begin{align*}
&( vf)(k_\beta ,k_\gamma )=
(2\pi)^{-\frac{3}{2}}{\int\limits_{({\T}^3)^2} }
 v (k_\beta-k_\beta' )
 \delta (k_\beta
+k_\gamma -k_\beta '-k_\gamma ')f(k_\beta ',k_\gamma')dk_\beta
'dk_\gamma',\\
  &f\in L^{(s)}_2(({\T}^3)^2),
  \end{align*}
where $\delta (\cdot)$ denotes the  Dirac delta-function at the
origin.

The functions $\varepsilon(k) \quad \mbox{and} \quad  v(k),\,
$  are  given by the Fourier series

\begin{equation*}
  \varepsilon
(k)=\sum_{j =1}^3(1-cosk^{(j)}), \quad  v
(k)=(2\pi )^{-3/2}\sum_{s\in
{{\Z}}^3} {\hat v}\,(s)\,e^{\mathrm{i}(k,s)},\,\,
\end{equation*}
 with
\begin{equation*}
(k,s)={\sum }_{j=1}^3 k^{(j)}s^{(j)}, \quad
k=(k^{(1)},k^{(2)},k^{(3)})\in {\T}^3,\quad
s=(s^{(1)},s^{(2)},s^{(3)})\in {{\Z}}^3.
\end{equation*}

The three-particle Hamiltonian $H$ is of  the form
\begin{equation*}
 H= H_0-{V}_{1}-V_{2}-V_{3},
\end{equation*}
where $H_0$ is the multiplication operator by the function
$\sum_{\alpha =1}^3 \varepsilon (k_\alpha)$
\begin{equation*}
( H_0f)(k_1,k_2,k_3)  = \sum_{\alpha =1}^3 \varepsilon
(k_\alpha)f(k_1,k_2,k_3),
\end{equation*}
\begin{align*}
&( Vf)(k_1,k_2,k_3)= (2\pi)^{-\frac{3}{2}}{\int\limits_{({\T}^3)^3}
} v (k_\beta-k_\beta')\times\\&\times \delta (k_\alpha -k_\alpha')\delta (k_\beta
+k_\gamma -k_\beta '-k_\gamma ')f(k'_1,k'_2,k'_3) dk'_1
dk_2'dk'_3, \quad f \in L^s_2(({\T}^3)^3).
\end{align*}

 \section{Decomposition of Hamiltonians into von Neumann direct integrals.
 Quasimomentum and coordinate systems}

 Given $m\in \N$, denote by $\widehat U^m_t$, $t\in {\Z}^3$ the unitary operators
 on the Hilbert space
  $\ell_2(({\Z}^3)^m)$ defined by:
\begin{equation*}
(\widehat U^m_tf)(n_1,n_2,..., n_m)=f(n_1+t,n_2+t,...,n_m+t),\quad f\in
\ell_2(({\Z}^3)^m).
\end{equation*}
We  easily see that
\begin{equation*}
  \widehat U^m_{t+\tau}=\widehat U^m_t \widehat U^m_ \tau,\quad  t, \tau \in \Z^3,
  \end{equation*}
i.e., $ \widehat U^m_t,t\in\Z^3 $ is a unitary  representation of
the abelian group $\Z^3$ in the Hilbert space
  $\ell_2(({\Z}^3)^m).$
 Since $\ell^{(s)}_2((\Z^{3})^{m})$ is invariant with respect to the group
  $\widehat U^m_t,t\in\Z^3,$
i.e.,  $$\widehat U^m_t  \ell^{(s)}_2((\Z^{3})^{m}) \subset
\ell^{(s)} _{2}((\Z^{3})^{m}),$$ we denote by   $\widehat U_{st}^m$
the restriction of $ \widehat  U^m_{t} $  to   the subspace $\ell
^{(s)}_2((\Z^{3})^{m})$.

Via the Fourier transform $ \cF_m^s$ the unitary representation of
$\Z^3$ in   $\ell^{(s)}_2(({\Z}^3)^m)$ induces a representation of
the group $\Z^3$ in the Hilbert space $L^{(s)}_2(({\T}^3)^m)$  by
unitary (multiplication) operators $  U_{st}^m=
(\cF_m^s)^{-1}\widehat U_{st}^m\cF_m^s$, $t\in \Z^3$ given by:
\begin{equation*}\label{grup}
(U_{st}^mf)(k_1,k_2,...,k_m)= \exp \big (-i(t,k_1+k_2+...+k_m)\big
)f(k_1,k_2,...,k_m),
\end{equation*}
\begin{equation*}
f\in L^{(s)}_2(({\T}^3)^m).
\end{equation*}

Denote by $K=k_1+k_2+...+k_m\in \T^3$ the {\it total quasi-momentum}
of the $m$ particles and define $\F_K^m$ as follows
\begin{align*}
&\F_K^m=\{(k_1,...,k_{m-1},K-k_1-...-k_{m-1}){\in
}({\T}^3)^m:\\
&k_1,k_2,...,k_{m-1}\in\T^3, K-k_1-...k_{m-1}\in\T^3\}.\end{align*}

Decomposing the Hilbert space $ L^{(s)}_2(({\T}^3)^m)$  into the
direct integral
\begin{equation*}
L^{(s)}_2(({\T}^3)^m) = \int_{K\in {\T}^3} \oplus
L^{(s)}_2(\F_K^m)dK
\end{equation*}
yields the  decomposition of the unitary representation $U_{st}^m,
\,t \in \Z^3$ into the  direct integral
\begin{equation*}
 U_{st}^m= \int_{K\in {\T}^3} \oplus  U_t(K)d K,
\end{equation*}
where
\begin{equation*}
U_t(K)=\exp(-i(t,K)) I \quad \text{on}\quad L^{(s)}_2(\F_K^m)
\end{equation*}
and $ I= I_{L^{(s)}_2(\F_K^m)}$ denotes the identity operator on the
Hilbert space $ L^{(s)}_2(\F_K^m)$.

The above Hamiltonians  $\hat h $ and $\widehat H$ obviously
commute with the groups of translations $ \widehat U_{st}^2,\, t\in
\Z^3$ and $ \widehat U_{st}^3, t\in \Z^d$ respectively, i.e.,
\begin{equation*}
 \widehat U_{st}^2 \hat h=\hat h \widehat U_{st}^2, \quad t\in \Z^3
\end{equation*}
and
\begin{equation*}
 \widehat U_{st}^3 \widehat H=\widehat  H \widehat U_{st}^3, \quad t\in \Z^3.
\end{equation*}
  Hence, the operators $ h$ and
  $ H$  can be decomposed into the direct integrals
\begin{equation}\label{fiber}
 h= \int\limits_{k \in {\T}^3}\oplus\tilde
h(k)d k\quad
\mbox{and}\quad
H=\int\limits_ {K \in {\T}^3}\oplus \widetilde H(K)dK
 \end{equation}
with respect to  the decompositions
\begin{equation*}L^{(s)}_2(\mathbb{T}^{3}) =
\int\limits_ {k \in {\T}^3} {\ \oplus } L^{(s)}_2(\F_K^2) d k,\quad \text{and}\quad
L^{(s)}_2((\mathbb{T}^{3})^{2}) = \int\limits_ {K \in {\T}^3} {\ \oplus }
L^{(s)}_2(\F_K^3) dK
\end{equation*}
respectively.

We introduce the mapping
\begin{equation*}
\pi^{(2)}:(\T^3)^2\to \T^3,\quad \pi^{(2)}((k_\beta,
k_\gamma))=k_\beta
\end{equation*}
resp.
\begin{equation*}
\pi^{(3)}:(\T^3)^3\to (\T^3)^2,\quad
\pi^{(3)}((k_\alpha, k_\beta, k_\gamma))=(k_\alpha,
k_\beta).
\end{equation*}

Denote by   $\pi^{(2)}_k$, $k\in
\T^3$  resp. $\pi^{(3)}_{K}$ , $K\in \T^3$   the restriction of $\pi^{(2)}$ resp.
$\pi^{(3)}$ onto $\F_k^2\subset (\T^3)^2$ resp.
  $\F_K^3\subset (\T^3)^3,$  i.e.,
\begin{equation*}\label{project}\pi^{(2)}_{k}= \pi^{(2)}\vert_{\F_k^2}\quad
\text{resp.}\quad \pi^{(3)}_{K}=\pi^{(3)}\vert_{\F_K^3}.
\end{equation*}
At this point it is useful to remark that $ \F^2_{k},\,\, k \in
{\bbT}^3 $ resp.
 $ \F^3_{K},\,\, K \in
{\bbT}^3 $ is three resp. six-dimensional manifolds isomorphic to
${\bbT}^3$ resp.
 ${({\bbT}^3)^2}.$

\begin{lemma}
The mapping   $\pi^{(2)}_{k}$, $k\in \T^3$ resp. $\pi^{(3)}_{K}$ ,
$K\in \T^3$ are bijective from $\F_k^2\subset (\T^3)^2$ resp.
$\F_K^3\subset (\T^3)^3$   onto $(\T^3)^2$ resp. $\T^3$ with the
inverse mapping given by
\begin{equation*}
(\pi^{(2)}_{k})^{-1}(k_\beta)=(k_{\beta},k-k_{\beta})
\end{equation*}
resp.
\begin{equation*}
(\pi^{(3)}_{K})^{-1}(k_\alpha,k_\beta)= (k_\alpha, k_{\beta},
K-k_\alpha-k_{\beta}).
\end{equation*}
\end{lemma}
\subsection{The fiber operators.}
The fiber operators $\tilde h(k),$ $k \in {\T}^3,$ from the direct
integral decomposition \eqref{fiber} are unitarily equivalent to the
operators $h(k),$ $k \in {\T}^3,$
 of the form
\begin{equation}\label{two} h(k) =h^{0}(k)-v.
\end{equation}

Let $L^{e}_2 ( {\T}^3)\subset L_2 ( {\T}^3)$ be the subspace of even functions.
The operators $h^{0}(k)$ and $v$ are defined on the Hilbert space $L^{e}_2 ( {\T}^3)$ by
\begin{equation*}(h^{0}
(k)f)(k_\beta)=\cE_{k}(k_\beta)f(k_\beta),\quad f \in L^{e}_2 ( {\T}^3),
\end{equation*}
where
\begin{equation}\label{E-alpha}
\cE_{k}(k_\beta)= \varepsilon (\frac{k}{2}-k_\beta) +\varepsilon
 (\frac{k}{2}+k_\beta)
\end{equation}
and
\begin{equation*} (v
f)(k_\beta)=(2\pi)^{-\frac{3}{2}} \int\limits_{{\T}^3}
v(k_\beta-k'_\beta) f(k_\beta')d k_\beta', \quad f \in L^{e}_2 ( {\T}^3).
\end{equation*}

The fiber operators $\widetilde H(K),$ \,$K \in {\T}^3$ from the
direct integral decomposition \eqref{fiber} are unitarily
equivalent to the operators $H(K),$ $K \in {\T}^3,$ and  are given by
\begin{equation*}
H(K)=H_0(K)-V_1-V_2- V_3.
\end{equation*}

 The operators $H_0(K)$
and $V_\alpha\equiv V ,\, \alpha=1,2,3,$ are defined on the Hilbert space $L^{e}_2 (({\T}^3 )^2)\cong L_2({\T}^3)\otimes L^{e}_2 ({\T}^3 )$
and  in the coordinates $(k_\alpha,k_\beta)\in({\T}^3)^2$ have form
\begin{equation}\label{TotalK}(H_0(K)f)(k_\alpha,k_\beta)=E
(K;k_\alpha,k_\beta)f(k_\alpha,k_\beta),\quad f\in  L^{e}_2 (({\T}^3 )^2) ,
\end{equation}
\begin{equation*} E(K;k_\alpha,k_\beta)= \varepsilon
(K-k_\alpha)+\varepsilon  (\frac{k_\alpha}{2}-k_{\beta}) + \varepsilon
(\frac{k_\alpha}{2}+k_{\beta})
\end{equation*}
and
 \begin{equation}\label{Poten}
V= I\otimes  v,
\end{equation}
wher $\otimes -$ is the tensor product and $ I=I_{L_2(\T^3)}$ is  the identity operator in $L_2(\T^3).$

 Later on all  our
calculations will be carried out in the configuration space, i.e.,
in a six-dimensional manifold $\F^3_{K}\subset (\T^{3})^{3}$
isomorphic to $(\T^{3})^{2}.$ As coordinates in $\F^3_{K}$ we will
choose one of the three pairs of vectors $(k_\alpha,k_\beta)$ which run independently
through the whole space $\T^{3}$ (if it does not lead to any
confusion we will write $(p,q)$ instead of $(k_\alpha,k_\beta)$.

 \section{The concept of a  zero energy resonance }
We  introduce the concept of a {\it zero energy resonance} ({\it
threshold resonance}) for the (lattice) two-particle operator $h(0)$
defined by \ref{two}.

We remark that under Hypothesis \ref{hypo} the sequence $\{\hat
v(s)\}_{s\in \Z^3}$ of the  Fourier coefficients  of the continuous
function
 $v(p)$ is an element of $\ell_2(\Z^3)$ and then
the equality

\begin{equation}\label{l2l2}
v(p)=(2\pi)^{-\frac{3}{2}}\sum_{s\in \Z^3}\hat
v(s)e^{\mathrm{i}(p,s)},
\end{equation} should be understood as follows:
the $L_2(\T^3)$-function $(2\pi)^{-\frac{3}{2}}\sum_{s\in \Z^3}\hat
v(s)e^{\mathrm{i}(p,s)}$ has the  continuous representative $v(p)$.

Since the function $\cE_k(q)$ has a unique non-degenerate minimum at
the point $q=\frac{k}{2}$ for any $k\in U_\delta (0)$ and $z\leq
\cE_{\text{min}}(k)$ the integral
\begin{equation}\label{yadr}
G(p,q;k,z)=\frac{1}{(2\pi)^{3}}\int\limits_{\T^3}
\frac{v^\frac{1}{2}(p-t)
 v^\frac{1}{2}(t-q)dt}
{\cE_k(t)-z}
\end{equation}
is finite, where
\begin{equation}\label{root}
v^\frac{1}{2}(p)=(2\pi)^{-\frac{3}{2}}\sum_{s\in \Z^3}\hat
v^\frac{1}{2}(s)e^{\mathrm{i}(p,s)}
\end{equation} is the kernel of
the operator $ v^\frac{1}{2}.$

We define the integral operator $G(k,z)$  in $L_2^e(\T^3)$ by
\begin{equation}\label{Bir-Sch}
G(k,z)f(p)=\int\limits_{\T^3} G(p,q;k,z)f(q)dq.
\end{equation}

\begin{lemma}\label{G.b-s} The operator $G(k,z),\,k\in U_\delta(0),$
$z \leq \cE_{\min}(k)$ acts in $ L_2^e(\T^3)$, is positive, belongs to
the trace class $\Sigma_1$, is continuous in $z$ from the left up to
$z=\cE_{\min}(k).$ \end{lemma}
\begin{proof} Lemma  \ref{G.b-s} can be proven in the same way as
Theorem 4.5
 in \cite{ALzM07}.\end{proof}
\begin{remark}\label{virt}
 Clearly, the operator $h(0)$ has an eigenvalue $z\le
\cE_0(0)=0$, i.e.,\\   $\Ker (h(0)-z I)\ne 0$, if and only if the
compact operator $G(0,z)$ on $L_2^e(\T^3)$ has an eigenvalue $1$ and
there exists a function $\psi\in \Ker(I-G(0,z))$ such that the
function $f$ given by
\begin{equation*}\label{fl2}
f(p)=\frac{(v^\frac{1}{2}\psi)(p)}{\cE_0(p)-z}
 \quad \text{a.e.}\quad p\in \T^3,
\end{equation*}
belongs to $L^2(\T^3)$. In this case $f\in \Ker(h(0)-z I)$.

Moreover, if $z<0$, then
\begin{equation*}\label{Kern}
\dim\Ker (h(0)-z I)=\dim \Ker (I-G(0,z))
\end{equation*}
and
$$
\Ker (h(0)-z I)=\{f\,\vert\,
f(\cdot)=\frac{(v^\frac{1}{2}\psi)(\cdot)}{\cE_0(\cdot)-z}\, , \,
\psi \in \Ker (I-G(0,z))\}.
$$
In the case of a threshold eigenvalue $z=0$ equality \eqref{Kern}
may fail to hold.

Therefore  equality \eqref{Kern} should be replaced by the
inequality
\begin{equation*}
\dim\Ker (h(0))\le\dim \Ker (I-G(0,0)).
\end{equation*}
\end{remark}
\begin{definition}\label{def}The operator $h(0)$ is said to have a zero energy resonance (at the
threshold) if $1$ is  an eigenvalue(single or multiple) of
$G(0,0)$ and at least one (up to a normalization) of the
associated eigenfunctions $\psi$ satisfies the condition
\begin{equation*}
\frac{(v^\frac{1}{2}\psi)(\cdot)}{\cE_0(\cdot)}\notin L^e_2(\T^3),
\end{equation*}
i.e.,
\begin{equation*}
1\le\dim \Ker (I-G(0,z))\ge \dim\Ker (h(0)-z I)+1.
\end{equation*}
\end{definition}

 \begin{remark}\label{zero-chala}
  The operator $h(0)$ has   a zero energy resonance, if and only if the point $1$ is a simple eigenvalue of
the operator $G(0,0)$ and the corresponding eigenfunction $\psi \in
L_2^e(\T^3)$ satisfies the condition $\int
v^{1/2}(p)\psi(p)dp\not=0. $
\end{remark}

\section{Statement of the main results}

We set:
\begin{align*}
E _{\min }(K)\equiv\min_{k_\alpha,k_\beta\in \T^3}E(K,k_\alpha,k_\beta),\quad E_{max
}(K)\equiv\max_{k_\alpha,k_\beta\in \T^3}E(K,k_\alpha,k_\beta).
\end{align*}
Our main results are as follows.

\begin{theorem}\label{mavjud}

Assume Hypotheses \ref{hypo}. Let  $h(0)$ have  a zero energy
resonance.
 Then  for all  $k \in U^0_\delta(0)$ the operator $h(k)$
has a unique positive eigenvalue $z(k)$ lying below the essential
spectrum. Moreovere $z(k)$ is analytic in $U^0_\delta(0).$
\end{theorem}

\begin{theorem}\label{ess}  For the essential spectrum
${\sigma }_{ ess } ( H(K) )$ of  $ H (K), K \in \T^{3} $ the following equality
holds
$$
\sigma_{ess}(H (K))=\bigcup _{p\in \T^3}\{
\sigma _d(h(p)) +\varepsilon(K-p)\}\cup [E _{min}(K),E _{max}(K)],
$$
where $\sigma _d(h(k))$ is the discrete spectrum of the operator
$h(k),  \,k\in \T^{3}$.
\end{theorem}

We denote by  $N(K,z)$ the number of eigenvalues of
$H(K),\,K\in \T^3$ below $z \leq \tau(K),$ where
\begin{equation}\label{essinf}
\tau(K)=\inf\sigma_{ess}(H (K)).
\end{equation}

\begin{theorem}\label{asimpZ} Assume Hypotheses \ref{hypo} and that the
 operator  $ h(0)$ has a zero energy resonance.

 Then the operator $H(0)$ has infinitely many eigenvalues lying
below the bottom of the essential spectrum and the function $N(0,z)$
obeys the relation
\begin{equation}\label{asimz} \lim\limits_{z \to
-0}\frac{N(0,z)}{|\log |z||}=\frac{\lambda _0}{2\pi},
\end{equation}
where   $\lambda _0$ the  unique positive solution of the equation
\begin{equation}\label{lamb}
 \lambda = \frac{8 \sinh \pi\lambda /6}{\sqrt 3 \cosh \pi\lambda/2}.
 \end{equation}

\end{theorem}

\begin{theorem}\label{asimpK} Let the conditions of Theorem \ref{asimpZ}
fulfilled. Then for all $ K \in U_\delta^0(0)$ the number $N(K,0)$
is finite and the following asymptotics holds
\begin{equation}\label{asimK}
\lim\limits_{|K| \to 0}\frac{N(K,0)}{|\log |K||}=
\frac{\lambda _0}{\pi}.
\end{equation}
\end{theorem}

\section{ Spectral properties of the two-particle operator $h (k)$ }

By  Weyl's theorem the essential spectrum
   $\sigma_{\text{ess}}(h(k))$ of the
operator $h(k),k \in \T^3$ defined by \eqref{two} coincides with the
spectrum $ {\sigma}( h_0 (k) ) $ of the non-perturbed operator
$h_0(k).$ More specifically,
\begin{equation*}
 \sigma_{\text{ess}}(h(k))= [\cE_{\min}(k)
,\cE_{\max}(k)], \end{equation*}
 where
\begin{equation*}\label{min}
 \cE_{\min }(k)\equiv\min_{p\in
\T^3}\cE_k(p),\quad \cE_{\max}(k)\equiv\max_{p\in {\T}^3}\cE_{k}(p)
\end{equation*} and $\cE_{k}(p)$ is
defined by \eqref{E-alpha}.

We denote by
$r_0(k,z)$ resp.  $r(k,z)$ the resolvent of the operator  $h_0(k)$ resp.  $h(k).$

Recall  the operators $r_0(k,z)$ and $r(k,z)$   are connected by the
following relations
\begin{equation}\label{resrel}
r(k,z)=r_0(k,z)+r_0(k,z)v r(k,z)= r_0(k,z)+r(k,z)v r_0(k,z).
\end{equation}

Set
\begin{equation*}
w(k,z)=I+v^{\frac{1}{2}} r(k,z)v^{\frac{1}{2}}.
\end{equation*}
   The resolvent relation \eqref{resrel} yields
\begin{equation*}
w(k,z)=(I-v^{\frac{1}{2}} r_0(k,z)v^{\frac{1}{2}})^{-1}.
\end{equation*}

 For a bounded self-adjoint operator $A,$
we define $n(\lambda,A)$ as
\begin{equation}\label{n()}
n(\lambda,A)=sup\{ dim F: (Au,u) > \lambda,\, u\in F,\,||u||=1\}.
\end{equation} The number $n(\lambda,A)$ is equal the infinity if
$\lambda$ is in the essential spectrum of $A$ and if $n(\lambda,A)$
is finite, it is equal to the number of the eigenvalues of $A$
bigger than $\lambda$.

The following lemma is  the Birman-Schwinger principle for the
two-particle Schr\"{o}dinger operators on the lattice $\Z^3$.
\begin{lemma}\label{G.b-s} For any $z\leq \cE_{\min }(k)$ the
following equality holds
\begin{equation*}
n(-z,-h(k)) =n(1,G(k,z)),
\end{equation*}
where $G(k,z),z\leq \cE_{\min }(k)$ is defined by \eqref{Bir-Sch}.
\end{lemma}
\begin{proof} Lemma  \ref{G.b-s} can be proven in the same way as
Theorem 4.5
 in \cite{ALzM07}.
\end{proof}

{\bf Proof of Theorem \ref{mavjud}}.
  By the assumptions of Theorem \ref{mavjud}
 the equation
\begin{equation*}\label{G0}
 G(0,0) \psi=
v^\frac{1}{2}r_0(0,0)v^\frac{1}{2}\psi=\psi
\end{equation*}
 has a nonzero solution  $\psi\in L_2^{e}(\T^3).$
One has
\begin{equation*}
(\psi,\psi)= (r_0(0,0)v^\frac{1}{2}\psi,v^\frac{1}{2}\psi)=
\int\limits_{{\bf T}^3}
\frac{|(v^\frac{1}{2}\psi)(p)|^2}{\cE_0(p)}dp,
\end{equation*}
with the continuous function $(v^\frac{1}{2}\psi)(\cdot)$ on $\T^3.$
  Then
\begin{align*}
&(v^\frac{1}{2}r_0(k,\cE_{\min }(k))v^\frac{1}{2}\psi,\psi)\\
&= \int\limits_{{\bf T}^3} \frac{|(v^\frac{1}{2}\psi)(p)|^2
dp}{\cE_k(p)-\cE_{\min }(k)}= \int\limits_{{\bf T}^3}
\frac{|(v^\frac{1}{2}\psi)(p)|^2dp}{\sum_{i=1}^3
2(1-cosp_i)cos(k_i/2)}\\
&>\int\limits_{{\T}^3}
\frac{|(v^\frac{1}{2}\psi)(p)|^2dp}{\sum_{i=1}^3 2(1-cosp_i)}=
\int\limits_{{\ T}^3}
\frac{|(v^\frac{1}{2}\psi)(p)|^2}{\cE_0(p)}dp=(\psi,\psi).
\end{align*}

By definition of \eqref{n()} this means that
$n(1,v^\frac{1}{2}r_0(k,\cE_{\min }(k))v^\frac{1}{2})>1.$ By the
Birman-Schwinger principle one concludes  that $h(k),k\neq0$ has an
eigenvalue lying below $\cE_{\min }(k).$

One checks that the equality
\begin{align*}\label{change}
&(h(k)f,f)=\int\limits_{{\T}^3}\varepsilon(\frac{k}{2}+q)\mid
f(q)\mid^2 d q  + \int\limits_{{\T}^3}\varepsilon(\frac{k}{2}-q)\mid
f(q)\mid^2 d q\\&-\int\limits_{{\T}^3}\int\limits_{{\T}^3}
v(p-q)f(q)\overline{f(p)}d q d p,\,f\in L_2^e(\T^3)
\end{align*}
holds. Since the function $f$ is even,  making a change of variables
in the integrals in the r.h.s.of the latter equality, we have
\begin{equation*}
(h(k)f,f)= (h(0)g,g)> 0,
\end{equation*}
 where $g(q)=f(\frac{k}{2}-q).$ This means that
$h(k)> 0$ for any $k \in \T^3\setminus\{0\}.$

The operator $G(k,0),\,k\in U_\delta(0)$ is analytic in
$U_\delta(0)$ and hence by the theorem \cite {RSIV} the eigenvalue
$z(k)$ is unique and analytic in $k \in U_\delta(0).$ $\Box$

Denote by $G_1$ the operator with the kernel

\begin{equation}\label{G1}
G_1(p,p')=-\frac{1}{4\pi} v^{\frac{1}{2}}(p)v^{\frac{1}{2}}(p').
\end{equation}

The following asymptotics plays a  crucial role in the proof of the
main result (Theorem \ref{asimpK} and \ref{asimpZ})(see also
\cite{Ltmf92}).
\begin{lemma}\label{rask} Assume Hypotheses \ref{hypo}.

(i) For all $k \in U_\delta(0)$ the following decomposition holds
\begin{equation}\label{G=}
G(k,0)=G(0,0)+\frac{1}{2}|k|G_1+|k|^{2}G_2(k) ,
\end{equation}
where  the operator $G_2(k)$ is continuous in $k\in U_\delta(0).$

(ii)  For all  $z\leq 0$ the following decomposition holds
$$
G(0,z)=G(0,0)+  G_1
(-z)^{\frac{1}{2}}+(-z)^{\frac{1}{2}+\theta_1} \tilde G_2(z)
,\,\,\,\theta_1<\theta,
$$
where the operator $\tilde G_2(z)$ is continuous in  $z\leq 0.$
\end{lemma}

\begin{proof}
For any $k \in U_\delta(0)$ the function  $\cE_k (q)$ has a unique
non-degenerate minimum at $q_k=\frac{k}{2}.$ Therefore,  by virtue
of the Morse lemma (see \cite{Fed}) there exists a regular
one-to-one mapping $q=\varphi(k;t)$ of a certain ball $W_\gamma(0)$
of radius $\gamma>0$ with the center at the origin to a neighborhood
$\tilde W(q_k)$ of the point $q_k$  such that
\begin{align}\label{eps=t2}
\cE_k(\varphi(k;t))=t^2+\cE_{\min }(k)
\end{align}
with $\varphi(k;0)=0$ and for the Jacobian $J(\varphi(k;t)) $ of the
mapping $q=\varphi(k;t)$ the equality
 \begin{align}\label{Jac}
J(\varphi(k,0))=\frac{1}{\sqrt{cos\frac{k^{(1)}}{2}cos\frac{k^{(2)}}{2}cos\frac{k^{(3)}}{2}}}
\end{align}
holds.

Since the function $\cE_0(\cdot)$ has a unique non-degenerate
minimum at $t=0$  by dominated convergence the finite limit
$$
G(p,q;0,0)=\lim_{k\to 0} G(p,q;k,0)
$$
exists for all $p,q\in \T^3.$

For all $p,q\in \T^3$  the following inequalities
\begin{align}\label{D-D}
&|G(p,q;k,0)-G(p,q;0,0)|\leq C |k|,
\end{align}
\begin{equation}\label{partial D}
\big |  \frac{\partial }{\partial |k|}G(p,q;k,0)- \frac{\partial}
{\partial |k|}G(p,q;0,0) \big |<C |k|^2,\quad k \in U_\delta(0)
\end{equation}
hold for some positive $C$ independent of $p$ and $q.$

Indeed, the function $G(p,q;\cdot,0)$  can be represented as
\begin{equation*}
G(p,q;k,0) =G_1(p,q;k,0)+G_2(p,q;k,0)
\end{equation*}
with
\begin{equation*}\label{I1}
G_1(p,q;k,0)=\frac{1}{(2\pi)^{3}}\int\limits_{\tilde W(q_k)}
\frac{v^\frac{1}{2}(p-t)
 v^\frac{1}{2}(t-q)dt}
{\cE_k (t)},\quad k \in U_\delta(0),
\end{equation*}
and
$$
G_2(p,q;k,0)=\frac{1}{(2\pi)^{3}}\int\limits_{{\T}^3\setminus \tilde
W(q_k)} \frac{v^\frac{1}{2}(p-t)
 v^\frac{1}{2}(t-q)dt}
{\cE_k (t)} \quad k \in U_\delta(0).
$$

Since for any $k \in U_\delta(0)$ the function $\cE_k (\cdot)$ is
continuous on the compact set\\ $\T^3 \setminus \tilde W(q_k)$ and
has a unique minimum at $q_k \in U_\delta(0)$ there exists
$M=const>0$ such that
$$\inf\limits_{q\in \T^3\setminus \tilde W(q_k)}\cE_k
(q)\geq M.$$

Then  for all $p,q\in \T^3$ we have
\begin{equation}\label{D2}
|G_2(p,q;k,0)-G_2(p,q;0,0)|\leq C k^2,\,k\in U_\delta(0)
\end{equation}
for some $C>0$ independent of $p,q\in \T^3.$

For all $p,q\in \T^3$ we consider
\begin{align}\label{D1}
G_1(p,q;k,0)-G_1(p,q;0,0)=-\frac{1}{(2\pi)^{3}} \int\limits_{\tilde
W(q_k)}
 \frac{(\cE_0(t)-\cE_k(t))v^\frac{1}{2}(p-t)
 v^\frac{1}{2}(t-q)dt}
{\cE_k(t)\cE_0(t)}.
\end{align}
In the integral in \eqref{D1} making a change of variable
$q=\varphi(k;t)$ and using  equality \eqref{eps=t2} we obtain
\begin{align}\label{D1-D1}
G_1(p,q;k,0)-G_1(p,q;0,0)=-\frac{\cE_{\min}(k)}{(2\pi)^{3}}
\int\limits_{ W_\gamma(0))} \frac{v^\frac{1}{2}(p-\varphi(k;t))
 v^\frac{1}{2}(\varphi(k;t)-q)J(\varphi(k;t))}
{t^2(t^2+\cE_{\min}(k))}dt.
\end{align}

Going over in the integral in \eqref{D1-D1} to spherical coordinates
$t=r\omega,$ we reduce it to the form
\begin{align*}\label{I2}
G_1(p,q;k,0)-G_1(p,q;0,0)=-\frac{\cE_{\min}(k)}{(2\pi)^{3}}
\int_0^{\gamma} \frac{F(p,q;k,r)}{r^2+\cE_{\min}(k)} dr,
\end{align*}
with
$$
F(p,q;k,r)=\int_{\Omega_2}v^\frac{1}{2}(p-\varphi(k,r\omega))
 v^\frac{1}{2}(\varphi(k,r\omega)-q)J(\varphi(k,r\omega))d\omega,
$$
where $\Omega_2$ is the unit sphere in $\R^3$ and $d \omega$ is the
element of the unit sphere in this space.  For all $p,q\in \T^3$ we
see that
\begin{equation}\label{GEld}
|F(p,q;k, r)-  F(p,q;0,0)|\leq C (r^{\theta} +|k|^{2})
\end{equation}
for some $C>0$ independent of $p,q\in \T^3.$

Indeed, applying equality \eqref{root} and taking into account that
the function $v^\frac{1}{2}(\cdot)$ is even on $\T^3$ we get

\begin{equation}\label{v-v}
|v^\frac{1}{2}(p-t)
 v^\frac{1}{2}(t-q)-v^\frac{1}{2}(p)
v^\frac{1}{2}(q)|\leq \frac{1}{(2\pi)^{3}}\sum_{s \in \Z^3}|\hat
v(s)||e^{-2\mathrm{i}(t,s)}-1|.
\end{equation}

For any $0<\theta \leq 1 $ the inequality
 $|e^{-2\mathrm{i}(t,s)}-1|\leq  C |t|^\theta |s|^\theta ,$
$p\in \T^3,\,\,s \in \Z^3$ holds. The Hypothesis \ref{hypo}
 and inequality \eqref{v-v} yield the inequality
\begin{equation*}\label{v-}
|v^\frac{1}{2}(p-t)
 v^\frac{1}{2}(t-q)-v^\frac{1}{2}(p)
v^\frac{1}{2}(q)|\leq    C |t|^\theta,\,\,\,\,\frac{1}{2}<\theta\leq
1,
\end{equation*}
for some  $C>0,$ independent on $p, q\in \T^3.$

Since the function  $\varphi(k,\cdot)$ (see. \cite{Fed}) is regular
and from  \eqref {Jac} we have
$$
 |J(\varphi(k,r))-J(\varphi(k,0))|\leq C |r|^2 \quad  \mbox{and} \quad
  |J(\varphi(k,0))-J(\varphi(0,0))|\leq C |k|^2.
$$
Thus
\begin{equation*}
|F(p,q;k, r)-  F(p,q;0,0)|\leq |F(p,q;k, r)-  F(p,q;k, 0)|+|F(p,q;k,
0)- F(p,q;0, 0)|\leq
\end{equation*}
\begin{equation*}
\int_{\Omega_2}|v^\frac{1}{2}(p-\varphi(k,r\omega))v^\frac{1}{2}(\varphi(k,r\omega)-q)-
v^\frac{1}{2}(p)v^\frac{1}{2}(q)||J(\varphi(k,0))|d\omega+
\end{equation*}
\begin{equation*}
\int_{\Omega_2}|v^\frac{1}{2}(p-\varphi(k,r\omega))v^\frac{1}{2}(\varphi(k,r\omega)-q)|
|J(\varphi(k,r\omega))-J(\varphi(k,0)|d\omega+
\end{equation*}
$$
\int_{\Omega_2}|v^\frac{1}{2}(p)||v^\frac{1}{2}(q)|
|J(\varphi(k,0))-J(\varphi(0,0))|d\omega\leq C (r^{\theta}
+|k|^{2}).
$$

For any $p,q\in \T^3$  the function $G_1(p,q;k,0)-G_1(p,q;0,0)$ can
be written in the form
\begin{equation}\label{F(0)}
G_1(p,q;k,0)-G_1(p,q;0,0)=- \frac{ F(p,q;0,0)}{(2\pi)^{3}}
\int_0^\gamma\frac{\cE_{\min}(k) dr}{r^2+\cE_{\min}(k)}-
\end{equation}
$$
\frac{1}{(2\pi)^{3}}\int_0^\gamma\frac{\cE_{\min}(k)(F(p,q;k,r)-F(p,q;0,0))}
{r^2+\cE_{\min}(k)}dr.
$$

By inequality \eqref{GEld} for all  $p,q\in \T^3$ we have
\begin{equation}\label{F-F}
\int_0^\gamma\frac{F(p,q;k,r)-F(p,q;0,0)}{r^2+\cE_{\min}(k)}dr\leq C
\int_0^\gamma\frac{r^\theta+|k|^2}{r^2+\cE_{\min}(k)}dr.
\end{equation}
The asymptotics
\begin{equation*}\label{m(k)}
\cE_{\min}(k)=\frac{1}{4}k^2+O(|k|^4) \quad \mbox{as} \quad  k\to 0
\end{equation*}
yields
\begin{equation*}\label{Limit1}
\frac{1}{k}\bigg (\int_0^\gamma\frac{\cE_{\min}(k)
dr}{r^2+\cE_{\min}(k)}- \frac{1}{4}\int_0^\gamma\frac{k^2
dr}{r^2+\frac{1}{4}k^2}\bigg )\to 0\,\,as\,\,k \to 0.
\end{equation*}
Computing the integrals
\begin{equation*}
\int_{0}^{\gamma}
  \frac{|k|}{r^2+\frac{1}{4} k^2}dr
\quad\mbox{and}\quad \int_{0}^{\gamma} \frac{
r^{\theta}+\frac{1}{4} |k|^2}{r^2+\frac{1}{4} k^2}dr
\end{equation*}
we obtain
\begin{equation*}\label{Limit}
\int_{0}^{\gamma}
  \frac{|k|}{r^2+\frac{1}{4} k^2}dr \to \pi \,\,as\,\,|k| \to 0
\quad\mbox{and}\quad \int_{0}^{\gamma}
  \frac{|k| (r^{\theta}+|k|^2)}{r^2+\frac{1}{4} k^2}dr\to 0 \,as\,\,|k| \to 0.
\end{equation*}

Using equality \eqref{F(0)} and inequality \eqref{F-F} we have
\begin{align*}\label{RhDer}
&\frac{\partial}{\partial |k|}G_1(p,q;0,0)=
 \lim_{|k|\to 0+}
  \frac{G_1(p,q;k,0)-G_1(p,q;0,0)}{|k|}=
 -\frac{1}{8\pi}  F(p,q;0,0).
\end{align*}
with the continuous function $F(p,q;0,0).$ Therefore the inequality
\begin{equation}\label{D1-baho}
\big |G_1(p,q;k,0)-G_1(p,q;0,0) \big |<C |k |,\quad k\in U_\delta(0)
\end{equation}
holds for some positive $C$ independent of $p,q\in \T^3.$

 Then from \eqref{D2} and \eqref{D1-baho}  it follows that the right-hand derivative of
 $G_1(p,q;\cdot, 0)$ at $|k|=0$ and
 for all $p,q\in\T^3$
we have
$$
 \frac{\partial }{\partial |k|}G(p,q;0,0)=
-\frac{1}{8\pi}v^\frac{1}{2}(p)v^\frac{1}{2}(q) .
$$

Comparing \eqref{D2} and \eqref{D1-baho} we get \eqref{D-D}. In  the
same way one can prove the inequality \eqref{partial D}. Then we
obtain \eqref{G=}.
 The part (ii) of Theorem \ref{rask} can be proven in  the same
 way as part (i)(see also \cite{Ltmf92}).

\end{proof}

\begin{theorem} \label{raswk}  Assume Hypotheses \ref{hypo}.

$a)$ If zero is a regular point of $h(0),$ then the operator
$w(k,0)$ resp. $w(0,z)$ is bounded for all $k \in U_\delta(0)$ resp.
$z \leq 0.$

$b)$ If $h(0)$ has a zero energy resonance, then the following
equalities hold

\begin{equation}\label{wraz}
w(k,0)= \frac{8\pi (\cdot,\psi)\psi}{\varphi^2(0)||\psi||^2 |k|}
 + w_1(k),\, k \in  U_\delta^{0}(0),\,\,
\end{equation}
and

\begin{equation}\label{w1/2}
w^{\frac{1}{2}}(k,0)= \frac{2\sqrt{2\pi}
(\cdot,\psi)\psi}{\varphi(0)||\psi||
 \sqrt {|k|} }
 +\tilde  w_1(k),\, k \in  U_\delta^{0}(0),
\end{equation}
resp.
\begin{equation*}\label{wrazz}
w(0,z)=\frac{4\pi (\cdot,\psi)\psi}{\varphi^2(0)||\psi||^2
\sqrt{-z}}
 +   w_2(z),\,\,z <0,
\end{equation*}
and
\begin{equation*}\label{w1/2raz}
w^{\frac{1}{2}}(0,z)= \frac{2\sqrt{\pi}
(\cdot,\psi)\psi}{\varphi(0)||\psi||
 \sqrt[4]{-z} }
 +\tilde   w_2(z),\,z<0,
\end{equation*}
where $\varphi(0)=(v^\frac{1}{2},\psi),\,$ the operators $ w_1(k)$
and $\tilde  w_1(k)$  are continuous in $k \in
 U_\delta(0)$ resp. $w_2(z)$ and  $\tilde w_2(z)$
 are continuous in $z \leq 0.$
\end{theorem}

\begin{proof} We shall prove the part of Theorem \ref{raswk}
concerning the operator $w(k,0).$  The part of Theorem \ref{raswk}
concerning  $w(0,z)$ can be proven in  the same way.

 $(a)$ Let zero
be a regular point. Since the operator $G(k,0)$ is analytic in
$U_\delta(0),$ compact and the number 1 is not an eigenvalue for
$G(0,0)$, the operator $(I-G(k,0))^{-1}$ exists for sufficiently
small $ |k|, k \in
 U_\delta(0).$

By the representation \eqref{G=} we have
$$
||w(k,0)||<C<\infty.
$$

$b).$
  Let the $h(0)$  have  a zero energy resonance.
   Denote by  $P_0$  the one dimensional projector onto the subspace
    ${\cH}_0$ associated with
$\psi$ and by $P_1$      the projector onto its orthogonal
complement ${\cH}_1$, so that $P_0 \oplus P_1=I.$ Let us write the
operator $A=A(k,0)=I-v^{\frac{1}{2}} r_0(k,0) v^{\frac{1}{2}}$ in
the matrix form:
$$
A=\left (
\begin{array}{ll}
A_{00} \,\,A_{01}\\
A_{10} \,\,A_{11}
\end{array}
\right ),
$$
where $A_{ij}=P_i A P_j: {\cH}_j \to {\cH}_i,\,\, i,j =0,1.$

Since $G(0,0)\psi=\psi$ , $P_0(I-G(0,0))P_1=0$, by Lemma \ref{rask}
 we have that
 $$
A_{00}= |k| \big ( -P_0 G_1P_0-|k|P_0 G_2(k)
 P_0 \big ),
$$
$$
A_{01}=|k| \big ( P_0 [-G_1- |k|G_2(k)]P_1 \big ),
$$
$$
A_{10}=A^*_{01},$$
\begin{equation}\label{A11}
A_{11}=
  P_1(I-G(0,0))P_1
-|k| P_1 [G_1+ |k| G_2(k)]P_1.
 \end{equation}

It is more convenient  to  consider the operator
$$
B=PAP,\quad P:= \left (
\begin{array}{ll}
\frac{P_0}{\sqrt{|k|}} &0\\
0  & P_1
\end{array}
\right ),
$$
instead of $A.$

Then from \eqref{A11} we obtain that
$$
B_{00}=-P_0 G_1P_0- |k|  P_0 G_2(k) P_0,
$$
$$
B_{01}=-\sqrt{|k|} P_0 [G_1+ |k| G_2(k)]P_1,
$$
$$
B_{10}=B^*_{01},
$$
$$
B_{11}= P_1(I-G(0,0))P_1 -|k| P_1 [G_1+|k| G_2(k)]P_1.
$$
Therefore $B=B^{(0)}+\tilde B,$ where
$$
B^{(0)}= \left (
\begin{array}{ll}
-P_0 G_1P_0 & 0\\
0 & P_1(I- G(0,0))P_1
\end{array}
\right )
$$
and $\tilde B=O(|k|^\frac{1}{2}) $ as $ k \to 0.$

 By the definition of $P_1$
the operator $F=(P_1(I-G(0,0))P_1)^{-1}$ exists in ${\cH}_1.$
 By
definition $P_0=||\psi||^{-2}(\cdot,\psi)\psi $, and we  obtain from
 $(v^{\frac{1}{2}},\psi)=\varphi(0)$ in
 \eqref{G1} that
$$
P_0G_1P_0=-\frac{1}{8\pi ||\psi||^2} P_0
(\psi,v^{\frac{1}{2}})^2=-\frac{\varphi^2(0)}{ 8 \pi ||\psi||^{2}}
P_0.
$$
Thus

$$(-P_0G_1P_0)^{-1}=\frac{8 \pi ||\psi||^{2}  }{ \varphi^2(0)} P_0=
\frac{8 \pi }{\varphi^2(0)}(\cdot,\psi)\psi.$$

 Now since
$B=(I+\tilde B(B^{(0)})^{-1})B^{(0)}$ and $\tilde B=O(|k|^\frac{1}{2})$ as $ k
\to 0$, we have
$$
B^{-1}=(B^{(0)})^{-1}+O(|k|^\frac{1}{2})= \left (
\begin{array}{ll}
\frac{8\pi(\cdot,\psi)\psi}{ \varphi^2(0)|k|}&0\\
0 &F
\end{array}
\right) +O(|k|^\frac{1}{2}). \textbf{}$$ Taking into account that
$w(k,0)=(A(k,0))^{-1}= PB^{-1}P$, we complete the proof of
\eqref{wraz}.

Let us prove \eqref{w1/2}. Since $r(k,0) \geq 0$  for $h(0) \geq 0$
we have $w(k,0) \geq I \geq 0$. Further, note that

$$
\left ( \frac{8\pi(\cdot,\psi)\psi}{ \varphi^2(0)||\psi||^2 |k|}
 \right )^{\frac{1}{2}}=
\frac{2\sqrt{2\pi}(\cdot,\psi)\psi}{ \varphi(0)||\psi||  \sqrt {|k|}
}
$$
and recall the well known inequality for arbitrary positive
operators A,\,B (see \cite{Bir}):
$||B^{\frac{1}{2}}-A^{\frac{1}{2}}||
 \leq
||B-A||^{\frac{1}{2}}.$ In combination with \eqref{wraz} this  yields
$$
\left | \left |  (w(k,0))^{\frac{1}{2}}-
\frac{2\sqrt{2\pi}(\cdot,\psi)\psi}{ \varphi(0)||\psi||
  \sqrt {|k|} }\right | \right |  \leq C .
$$

\end{proof}

\section{Spectrum of the operator $ H(K)$}
Since the particles are identical we have only one channel operator
 $H_{ch}(K), K {\ \in } {\T}^3$ acting in the
Hilbert space $L^{e}_2 (({\T}^3 )^2)\cong L_2({\T}^3)\otimes L^{e}_2
({\T}^3 ) $ as
\begin{equation*}
 H_{ch}(K)=H_0(K)- V,
\end{equation*}
where  $H_0(K)$  resp.  $V$ is defined by  \eqref{TotalK}  resp.  \eqref{Poten}.

The decomposition of the space $L^{e}_2 (({\T}^3 )^2)$ into the
direct integral

 $$L^{e}_2 (({\T}^3 )^2) = \int\limits_{k\in {\T}^3}
\oplus L^{e}_2({\T}^3) dk
$$
yields for the operator $H_{ch}(K)$ the decomposition into the
direct integral
 $$ H_{ch}(K) =  \int\limits_{k\in {\T}^3}
 \oplus H_{ch}(K,k) dk.$$
 The fiber operator $ H_{ch}(K,k) $ acts in the
Hilbert space $L^{e}_2({\T}^3)$ and has the form
\begin{equation*}
 H_{ch}(K,k) =h(k)+\varepsilon(K-k) I,
\end{equation*} where $ I$  is
identity operator and $h(k)$ is the two-particle operator defined
  by \ref{two}.
The representation of the operator $H_{ch}(K,k)$ implies the
equality
\begin{align*}\label{stucture}
 &\sigma (H_{ch}(K,k))
 = \sigma _d(h(k))\cup \big[\cE_{\min
}(k),\cE_{\max}(k) \big]+\varepsilon(K-k),
\end{align*}
 where $\sigma _d(h(k))$
 is the discrete spectrum of the operator
 $h(k).$
 \begin{lemma}\label{spec} The following
equality holds
\begin{align*}
&\sigma (H_{ch}(K))\\&=\cup _{k\in {\T}^3} \left \{ \sigma _d(h(k)+
\varepsilon(K-k) \right\} \cup [E_{\min}(K),E_{\max}(K)].
\end{align*}
\end{lemma}
\begin{proof}
The theorem (see, e.g.,\cite{RSIV}) on the spectrum of decomposable
operators and above obtained structure for the spectrum of $H_{ch}
(K,k)$  complete the proof.
\end{proof}

\begin{lemma}\label{equality} The following  equality holds
$$\sigma (H_{ch}(K))=\sigma _{ess}(H(K)).$$
\end{lemma}
\begin{proof} Lemma  \ref{equality} can be proven in  the same way as
Theorem 3.2
 in \cite{ALzM07}.
\end{proof}

{\bf Proof of Theorem \ref{ess}.} Lemma \ref{spec} and
\ref{equality} yield the proof of Theorem  \ref{ess}. $\Box$

\section{Birman-Schwinger principle for the operator $H(K).$}

Recall that (see \eqref{Poten})  the operator $V$ acting in
$L_2^{e}(({\T}^3)^2)$ has form
$$ Vf (k_\alpha,k_\beta)=(I \otimes
v) f (k_\alpha,k_\beta)$$

and hence the operator $V^{\frac{1}{2}}$ is of the form
$$ V^{\frac{1}{2}} f (k_\alpha,k_\beta)=(I \otimes  v ^{\frac{1}{2}})
f (k_\alpha,k_\beta)
$$

Let $W(K,z)$ and  $W^{\frac{1}{2}} (K,z), \,\,K \in \T^3, \,\, z \leq \tau (K),$
be the
operators in $L_2^{e}(({\T}^3)^2)$ defined as

\begin{equation*}\label{wrazz}
W(K,z)f(k_\alpha,k_\beta) =\big (I \otimes w(k_\alpha ,z-
\varepsilon( K-k_\alpha)) \big )
 f (k_\alpha,k_\beta),
\end{equation*}

\begin{equation}\label{Wraz}
W^{\frac{1}{2}} (K,z) f (k_\alpha,k_\beta) =\big (I \otimes
w^{\frac{1}{2}}
 (k_\alpha ,
z- \varepsilon( K-k_\alpha))  \big )
 f (k_\alpha,k_\beta).
\end{equation}

The operator $W(K,z), \,\,K \in \T^3, \,\, z \leq \tau (K),$
as related to  the resolvent
$R_{ch}(K,z)$ as

$$
W(K,z)=I+V^{\frac{1}{2}} R_{ch}(K,z)V^{\frac{1}{2}},
 $$
where  $R_{ch}(K,z)$ is  the resolvent of $H_{ch}(K), K \in \T^3.$

 One checks that
\begin{equation*}\label{decom}
 W(K,z)=(I-V^{\frac{1}{2}}
R_0(K,z)V^{\frac{1}{2}})^{-1}
\end{equation*}
where $R_0(K,z)$ the resolvent of the operator $H_0(K).$

For $z < \tau(K), \,\,K \in \T^3, $ the operator $W(K,z)$ is  positive.

Let  $
 {\bf T}(K,z), \,\,K \in \T^3, \,\, z \leq \tau (K),
$
be the operator in $L_2^{e}(({\T} ^3)^2)$  defined by

\begin{equation}\label{T(Fad)}
 {\bf T} ( K, z) = 2 W^{\frac{1}{2}}(K,z)V^{\frac{1}{2}}
R_0(K,z)V^{\frac{1}{2}}W^{\frac{1}{2}}(K,z), \,\,K \in \T^3, \,\, z \leq \tau (K).
\end{equation}

 By the definition of
$N(K,z), \,\,K \in \T^3, \,\, z \leq \tau (K),$ we have
$$
N(K,z)=n(-z,-H(K)),\,-z > -\tau(K).
$$
The following lemma is a realization of well known Birman-Schwinger
principle for the three-particle Schr\"{o}dinger operators on the
lattice $\Z^3 $(see \cite{Sob,Tam94} ).
\begin{lemma}\label{b-s}
The operator ${\bf T}(K,z), \,\,K \in \T^3, \,\, z < \tau (K),$
is compact and continuous in $z < \tau (K)$ and
the following equality holds
$$
N(K,z)=n(1,{\bf T}(K,z)).
$$
\end{lemma}
\begin{proof} We first verify the equality
\begin{equation}\label{tenglik}
N(K,z)=n(1,3R^{\frac{1}{2}}_0(K,z)VR^{\frac{1}{2}}_0(K,z)).
\end{equation}
 Assume that $u \in
{\cH}_{-H(K)}(-z)$, i.e., $((H_0(K)-z)u,u) < 3( Vu,u).$ Then
$$ (y,y) < 3(R^{\frac{1}{2}}_0(K,z) V
R^{\frac{1}{2}}_0(K,z)y,y),\quad y=(H_0(K)-z)^{\frac{1}{2}}u.
$$
Thus $N(K,z) \leq
n(1,3R^{\frac{1}{2}}_0(K,z)VR^{\frac{1}{2}}_0(K,z))$. Reversing the
argument we get the opposite inequality, which proves
\eqref{tenglik}. Any nonzero  eigenvalue of
$R^{\frac{1}{2}}_0(K,z)V^{\frac{1}{2}}$ is an eigenvalue for
$V^{\frac{1}{2}}R^{\frac{1}{2}}_0(K,z)$
 with  the same
algebraic and geometric multiplicities.

Therefore  we get
$$
n(1,3R^{\frac{1}{2}}_0(K,z) V R^{\frac{1}{2}}_0(K,z))=
n(1,3V^{\frac{1}{2}} R_0(K,z)V^{\frac{1}{2}}).
$$

Let us check that
$$
n(1,3R^{\frac{1}{2}}_0(K,z) V R^{\frac{1}{2}}_0(K,z))=n(1,{\bf
T}(K,z)).
$$
We shall show that for any $u \in {\cH}_{3R^{\frac{1}{2}}_0(K,z)V
R^{\frac{1}{2}}_0(K,z))}(1)$ there exists $y \in{\cH}_{{\bf
T}(K,z)}(1)$ such that $(y,y)<({\bf T}(K,z)y,y).$ Let $u \in
{\cH}_{3R^{\frac{1}{2}}_0(K,z)V R^{\frac{1}{2}}}(1)$ i.e.,
$$
(u,u)< 3(V^{\frac{1} {2}} R_0(K,z)V^{\frac{1}{2}} u,u)
$$
and hence
\begin{equation}\label{coordinate}
((I-V^{\frac{1}{2}} R_0(K,z) V^{\frac{1}{2} })u,u)< 2(V^{\frac{1}{2}
}R_0(K,z)V^{\frac{1}{2}} u,u).
\end{equation}
Setting  $y=(I- V^{\frac{1}{2}}R_0(K,z)V
^{\frac{1}{2}})^{\frac{1}{2}} u $ we have
$$
(y,y)< (2W^{\frac{1}{2}}(K,z)V^{\frac{1}{2}} R_0(K,z)
V^{\frac{1}{2}} W^{\frac{1}{2}} (K,z)y,y),
$$
i.e., $ (y,y)\leq ({\bf T}(K,z)y,y). $ Thus $
n(1,3R^{\frac{1}{2}}_0(K,z)V R^{\frac{1}{2}}_0(K,z)) \leq n(1,{\bf
T}(K,z)). $

In  the same way one checks  that  $ n(1,{\bf T}(K,z)) \leq
n(1,3R^{\frac{1}{2}}_0(K,z)V R^{\frac{1}{2}}_0(K,z))  . $

Finally we note that for any  $z<\tau(K)$ the operator $T(K,z)$ is
compact and continuous in $z.$

\end{proof}
\begin{remark}
On the left hand side of  \eqref{coordinate} the operator
$V^{\frac{1}{2}} R_0(K,z) V^{\frac{1}{2}}$ is a partial integral
operator, since the operator
$$V^{\frac{1}{2}} f
(k_\alpha,k_\beta)=V^{\frac{1}{2}}_\alpha f
(k_\alpha,k_\beta)=(I\otimes v^{\frac{1}{2}}) f (k_\alpha,k_\beta)
$$ is written in the coordinate $(k_\alpha,k_\beta)$ i.e., it is
an integral operator with respect to  $k_\beta.$

The right hand side of \eqref{coordinate} can be written as sum of
$V_1^{\frac{1}{2}} R_0(K,z) V_2^{\frac{1}{2}}$ and
$V_1^{\frac{1}{2}} R_0(K,z) V_3^{\frac{1}{2}},$ where the operator
$V=V_1$ is written in the coordinate $(k_\alpha,k_\beta)$,i.e., it
is integral operator  with respect to  $k_\beta.$ But the operators
$V=V_2$ and $V=V_3$ in the coordinates $(k_\beta,k_\alpha)$
constitute it is an integral operators with respect to  $k_\alpha$
and hence the operator $2V^{\frac{1}{2}} R_0(K,z) V^{\frac{1}{2}}$
on the  right hand side of \eqref{coordinate} is an  integral
operator.
\end{remark}

\section{The number of eigenvalues of the  operator
$H(K)$}

In this section we shall prove Theorems \ref{asimpZ} and
 \ref{asimpK}.
First we prove that the number $N(K,0)$ is finite.

\begin{theorem}\label{main} The following equality holds
$$
\lim\limits_{|K|\to 0} \frac{n(1,{\bf T}(K,0))} {|log|K||}
=\frac{\lambda _0}{\pi} \quad \mbox {resp.} \quad
\lim\limits_{|z|\to 0} \frac{n(1,{\bf T}(0,z))} {|log|z||}
=\frac{\lambda _0}{2\pi}
$$
where $\lambda _0$ is defined in Theorem \ref{asimpZ}.
\end{theorem}

Theorem \ref{main} will be deduced by a perturbation argument based
on Lemma 4.7, which has been proven in \cite{Sob}. For completeness,
we here reproduce the lemma.

\begin{lemma}\label{comp.pert}
Let $A (z)=A_0 (z)+A_1 (z),$ where $A_0(z)$ (resp.$A_1(z)$) is
compact and continuous in $z<0$ (resp.$z\leq 0$).  Assume that for
some function $f(\cdot),\,\, f(z)\to 0,\,\, z\to 0-$ one has
$$
\lim_{z\rightarrow 0-}f(z)n(\lambda,A_0 (z))=l(\lambda),
$$
and $l(\lambda)$ is continuous in $\lambda>0.$ Then the same limit
exists for $A(z)$ and
$$ \lim_{z\rightarrow 0-}f(z)n(\lambda,A (z))=l(\lambda).
$$
\end{lemma}

\begin{remark} According
to Lemma \ref{comp.pert} any perturbation of the operator $A_0(z)$
defined in Lemma \ref{comp.pert}, which is compact and continuous up
to $z=0$ does  not contribute to the asymptotics \eqref{asimz}.
%Throughout the proof of the following theorem we shall use this fact
%without further comments.
\end{remark}

Let the operator $h(0)$ have   a zero energy resonance ,i.e., the
number $1$ is an eigenvalue for $G(0,0)$ and $\psi$ is an associated
eigenfunction.

Let $\Psi:L_2({\T}^3)\rightarrow L_2^{e}(({\T}^3)^2)$ be the
operator defined by

\[ (\Psi f)(k_1,k_2)=
\frac{1}{||\psi||}{\psi(k_2)}f(k_1)
\]

and  $\Psi^* :L_2^{e}(({\T}^3)^2)\rightarrow L_2({\T}^3)$  be
adjoint operator,i.e.,
\[ ( \Psi^*
f)(k_1)=\frac{1}{||\psi||}\int\limits_{{\bf T}^3} \psi(k^{\prime}_2)
f(k_1,k^{\prime}_2)dk^{\prime}_2.
\]

Let $ \Delta^{-1/4}$ be the operator of multiplication by the
function $$(\cE_{\min}(p)+\varepsilon(K-p)-z)^{-1/4}.$$
  Theorem
\ref{raswk} and equality \eqref{Wraz} yield
\begin{equation}\label{W(K)}
W^{\frac{1}{2}}(K,0)=\frac{2\sqrt{\pi} ||\psi||}{ \varphi (0)} \Psi
\Delta^{-1/4} \Psi^{*} +\widetilde W_1(K),\,
\end{equation}
resp.
\begin{equation}\label{W(z)}
W^{\frac{1}{2}}(0,z)=\frac{2\sqrt{\pi} ||\psi||}{ \varphi (0)} \Psi
\Delta^{-1/4} \Psi^{*} +\widetilde W_2(z),\,
\end{equation}
 where $\widetilde W_1(K)$  resp. $\widetilde W_2(z)$ is a
continuous operator in $ K \in U_\delta(0) $ resp. in $z \leq 0.$

We define the operator $ {\bf T}^{(1)} (K,z)$ acting in $L_2^{e}(({\bf T}
^3)^2)$ by
\begin{equation}\label{Gamma}
{\bf T}^{(1)} (K,z)=\frac{8\pi ||\psi||^2}{\varphi^2 (0)}\Psi
\Delta^{-1/4}
 \Psi^*V^{\frac{1}{2}} R_0(K,z)
V^{\frac{1}{2}}\Psi   \Delta^{-1/4}
  \Psi^*.\end{equation}
\begin{lemma}\label{b}
 For all $ K \in U_\delta(0) $ resp.  $z \leq 0$ the operator
$ {\bf T}(K,0)-{\bf T}^{(1)} (K,0)$ resp. $ {\bf T}(0,z)-{\bf
T}^{(1)} (0,z)$
 is a compact operator.
\end{lemma}
\begin{proof}
 From \eqref{W(K)}, \eqref{Gamma}  and \eqref{T(Fad)} we have
$$
{\bf T}(K,0)-{\bf T}^{(1)}(K,0)=\frac{4\sqrt{\pi} ||\psi||}{\varphi
(0)}
 \Psi \Delta^{-\frac{1}{4}} \Psi^*V^{\frac{1}{2}}
R_0(K,0)V^{\frac{1}{2}}\,\,\, \widetilde W_1(K)+
$$
\begin{equation}\label{Gamma+}
+\frac{4\sqrt{\pi} ||\psi||}{\varphi (0)} \widetilde W_1(K)\,V^{\frac{1}{2}}
 R_0(K,0)V^{\frac{1}{2}}\Psi \Delta^{-\frac{1}{4}}
 \Psi^* +
2\widetilde W_1(K)\,V^{\frac{1}{2}}R_0(K,0)V^{\frac{1}{2}}\,\,\,
 \widetilde W_1(K).
\end{equation}
The kernels of the operators $\Psi \Delta^{-\frac{1}{4}}
\Psi^*V^{\frac{1}{2}} R_0(K,0)V^{\frac{1}{2}} $, $ V^{\frac{1}{2}}
R_0(K,0)V^{\frac{1}{2}}\Psi \Delta^{-\frac{1}{4}} \Psi^* $ and
$V^{\frac{1}{2}}R_0(K,0)V^{\frac{1}{2}}$ are bounded by
$$
\frac{C_1}{|p|(p^2+q^2)},\quad \frac{C_2}{(p^2+q^2)|q|}\quad and
\quad \frac{C_3}{(p^2+q^2)},
$$
respectively.

By passing on to a spherical  coordinates  system one checks that
the functions  $ |p|^{-1}(p^2+q^2)^{-1},\quad
(p^2+q^2)^{-1}|q|^{-1}\quad$ and $\quad (p^2+q^2)^{-1}$ are square
integrable.  By Lemma \ref{raswk} the operator $\widetilde W_1(K)$
is uniformly bounded in $K \in U_\delta(0) .$ Therefore the
right-hand side of \eqref{Gamma+} is  a compact operator for any
$ K
\in U_\delta(0) .$

In  the same way one checks that  $ {\bf T}(0,z)-{\bf
T}^{(1)}(0,z),\,$$z \leq 0,$ belongs to the Hilbert-Schmidt class.

 The proof of Lemma \ref{b-s} is completed.
\end{proof}
%%%%%%%%%%%%%%%%

For all $ K \in U_\delta(0) $ and $z \leq 0$  we define the operator
${T}^{(1)} (K,z):L_2({\T}^3)\rightarrow L_2({\T}^3),$ by
\begin{equation}\label{TPsi}
{T}^{(1)}(K,z)=\frac{8\pi ||\psi||^2}{\varphi^2(0)}\Delta^{-1/4}
\Psi^*V^{\frac{1}{2}} R_0(K,z) V^{\frac{1}{2}} \Psi \Delta^{-1/4}.
\end{equation}
 The kernel of  ${T}^{(1)}(K,z)$ in the coordinates $(p,q)$
  has form:
$$\frac{1}{\varphi^2(0)\pi^2}\frac{(E(K,\frac{2K}{3}+p,\frac{p}{2}+q)-z)^{-1}
 \varphi(\frac{p}{2}+q)\varphi(p+\frac{q}{2})}
 {
(\cE_{\min}(\frac{2K}{3}+p)+\varepsilon(\frac{K}{3}-p)-z)^{1/4}
(\cE_{\min}(\frac{2K}{3}+q)+\varepsilon(\frac{K}{3}-q-z)^{1/4}}$$

\begin{lemma}\label{T-cT}
  The discrete spectrum of ${\bf T}^{(1)}(K,0)$ and ${T}^{(1)}
(K,0)$  resp. ${\bf T}^{(1)}(0,z)$ and ${T}^{(1)} (0,z)$ coincides.
\end{lemma}
\begin{proof}
According to \eqref{TPsi} and \eqref{Gamma} we get
\begin{equation}\label{cT=}
 {\bf T}^{(1)}(K,0)=\Psi T^{(1)}(K,0)\Psi^*.
\end{equation}
Since  any nonzero eigenvalue of $\Psi^* T^{(1)}(K,0)\Psi$ resp.
$\Psi^* T^{(1)}(0,z)\Psi$  is an eigenvalue of
$\Psi\Psi^*T^{(1)}(K,0)$  as well,  with the same algebraic and
geometric multiplicities,  and $\Psi\Psi^*=I,$ where $I$ is the
identity  operator on $L_2(\T^3),$  we have $ n(1,{\bf T}^{(1)}
(K,0))=n(1, {T}^{(1)}(K,0) ) .$

On the other hand $\Psi^*\Psi $ is the identity operator in $
L_2(\T^3)$ and hence using  equality \eqref{cT=} we have
$ \sigma_d({\bf T}^{(1)} (K,0))=\sigma _d(T^{(1)}(K,0)), $

In  the same way one proves that  $ \sigma _d({\bf T}^{(1)}
(0,z))=\sigma _d(T^{(1)}(0,z)) .$ The Lemma \ref{T-cT} is proven.
\end{proof}

 \begin{lemma}\label{EK.det}
The function $E(K,\frac{2K}{3}+p,\frac{p}{2}+q)$ resp.
$\cE_{\min}(\frac{2K}{3}+p)+\varepsilon(\frac{K}{3}-p)$ has the
following asymptotics

 \begin{align}\label{eps exp}
E(K,\frac{2K}{3}+p,\frac{p}{2}+q)=\frac{K^2}{6}+p^2+(p,q)+q^2+O(|p|^4)+
 O(|q|^4)+O(|K|^4)
 \to
0\end{align}
as  $K, p,q\rightarrow 0$

resp.
\begin{align}\label{det exp}
\cE_{\min}(\frac{2K}{3}+p)+\varepsilon(\frac{K}{3}-p)=
\frac{K^2}{6}+\frac{3p^2}{4}+O(|K|^4)+O(|p|^4),
\end{align}
as  $K, p\rightarrow 0.$

\end{lemma}

\begin{proof}
 The asymptotics
\begin{equation}\label{eps. exp}
  \varepsilon(p)=\frac{1}{2}p^{2}+O(|p|^{4}) \quad as \quad
  p\rightarrow 0
\end{equation}
of the function $\varepsilon(p)$ yields \eqref{eps exp}.
 The
definition of $\cE_{\min}(k)$ and the representation \eqref{E-alpha}
gives the asymptotics
\begin{equation}\label{mkraz}
 \cE_{\min}(k)=\frac{1}{4}k^{2} +O(|k|^{4}) \quad
  as \quad k\rightarrow 0,
\end{equation}
which yields  \eqref{det exp}.
\end{proof}

Denote by $ \chi_\delta(\cdot)$ the  characteristic function of $
U_\delta(0)=\{ p\in \T^3:\,\, |p|<\delta \}.$

For all $ K \in U_\delta(0) $ and $z \leq 0$  we define the  operator
${T}^{(1)}(\delta,\frac{K^2}{6}+|z|),$
on $L_2({\T}^3)$ with the kernel
\begin{align*}
\frac{1}{\pi^2}
  \frac{ \chi (p)\chi (q)}
{ (\frac{3
p^2}{4}+\frac{K^2}{6}+|z|)^{1/4}(p^2+(p,q)+q^2+\frac{K^2}{6}+|z|)(\frac{3 q^2}{4}+\frac{K^2}{6}+|z|)^{1/4}
 }.
\end{align*}

\begin{lemma} \label{1Raznost}  The operator $ {T}^{(1)}(K,0)-{T}^{(1)} (\delta;
\frac{K^2}{6})$ resp. $ {T}^{(1)}(0,z)-{T}^{(1)}(\delta;|z|)$
belongs to the Hilbert-Schmidt class and is continuous in $K\in
U_\delta (0)$ resp. $z\leq 0.$
\end{lemma}
\begin{proof}

Applying the  asymptotics \eqref{eps exp} and \eqref{det exp}
one can estimate
the kernel of the operator $T (K,0) -T (\delta; \frac{K^2}{6})$
 by
\begin{equation*}
 C [ (p^2+q^2)^{-1} +
|p|^{-\frac{1}{2}}(p^2+q^2)^{-1} +
(|q|^{-\frac{1}{2}}(p^2+q^2)^{-1} + \frac{
|p|^{\theta}+|q|^\theta } {
|q|^{\frac{1}{2}}|p|^{\frac{1}{2}}(p^{2}+q^2 )}+1 ]
\end{equation*}
 and hence
the operator ${T}^{(1)}(K,0)-{T}^{(1)} (\delta; \frac{K^2}{6})$
 belongs to the Hilbert-Schmidt class for all $K\in U_{\delta}(0).$
 In combination with the continuity of the kernel of
the operator in $K\in U_\delta (0)$  this  gives   the continuity of
${T}^{(1)} (K,0)-{T}^{(1)}(\delta;\frac{K^2}{6})$ in $K\in U_\delta
(0).$

In  the same way one checks that  $
{T}^{(1)}(0,z)-{T}^{(1)}(\delta;|z|)$ belongs to the Hilbert-Schmidt
class and is continuous in  $z\leq 0.$

\end{proof}

Let
\begin{align*}
&{\bf S}({\bf r}):L_2((0,{\bf r}),
{\sigma_0})\to L_2((0,{\bf r}),{\sigma_0}),\,\,
{\sigma_0}=L_2(\S^2),
\\&{\bf r}=1/2 | \log \frac{K^2}{6}|\quad
\mbox{resp.}\quad {\bf r}=1/2 | \log |z||,\\
\end{align*}
$\S^2-$ being the unit sphere in $\R^3$, be the integral operator
with the kernel
\begin{align*}\label{Sobolov}
&  S(t;y)=\frac{2}{\sqrt 3 \pi^2 }\frac{1}{\cos h y+\frac{1}{2 } t},\\
 &y=x-x',\,x,x'\in (0,{\bf r}),\quad t=<\xi,
\eta>,\,\xi, \eta \in \S^2,\no
\end{align*}
and let
$$\hat{\bf S}(\lambda):\,\,\sigma_ 0\rightarrow
\sigma_0,\,\,\,\lambda\in (-\infty,+\infty) $$ be the integral
operator with the following kernel

\begin{equation}\label{Stlam}
\hat { S}(t;\lambda)=\int\limits_{-\infty}^{+\infty}\exp{\{-i\lambda
r\}}{S}(t;r)dr=\frac{1}{\sqrt 3 \pi} \frac{\sinh[\lambda(arc\cos
\frac{1}{2}t)]} {(1-\frac{1}{4}t^2)^{\frac{1}{2}}\sinh (\pi\lambda)}
\end{equation}

 For $\mu>0,$ define
\begin{equation}\label{sobU}
 {U}(\mu)= (4\pi)^{-1}
\int\limits_{-\infty}^{+\infty} n(\mu,\hat{\bf S}(y))dy.
\end{equation}
\begin{lemma}\label{sobol}
 The function  $U(\mu)$ is continuous in $\mu>0$,
 the following limit $$ \lim\limits_{{\bf r}\to \infty} \frac{1}{2}{\bf
r}^{-1}n(\mu,{\bf S}({\bf r}))={U}(\mu)$$ exists.
\end{lemma}

\begin{remark}
 This lemma can be proven in  the same way as
  the corresponding results of \cite{Sob}. In particular, the
continuity of ${U}(\mu)$ in $\mu>0$ is a result of Lemma 3.2,
Theorem 4.5  states  the existence of the limit
$$ \lim\limits_{{\bf r}\to
\infty} \frac{1}{2}{\bf r}^{-1}n(\mu,{\bf S}({\bf
r}))={U}(\mu).$$
\end{remark}

\begin{lemma}\label{main1} The equalities
$$
\lim\limits_{|K|\to 0} \frac{n(1,{T}^{(1)}(\delta,\frac{K^2}{6}))} {|log|K||}
=\frac{\lambda _0}{\pi}
$$
and
$$
\lim\limits_{|z|\to 0} \frac{n(1,{T}^{(1)}(\delta,|z|)} {|log|z||}
=\frac{\lambda _0}{2\pi}
$$
hold, where $\lambda_0$ is unique positive solution of the equation
\eqref{lamb}.
\end{lemma}

\begin{proof}
 The space of functions  having support in $U_\delta(0)$
 is an invariant subspace
for the operator $T^{(1)}(\delta,\frac{K^2}{6}).$

Let $T^{(1)}_{res}(\delta ,\frac{K^{2}}{6}),$  be the restriction of
the operator $T^{(1)}(\delta ,\frac{K^{2}}{6})$  on the invariant
subspace $L_{2}(U_{\delta }(0)).$

The operator $T^{(1)}_{res}(\delta ,\frac{K^{2}}{6})$  is unitarily
equivalent with the operator
 $T^{(2)}(\delta ,\frac{K^{2}}{6})$
  acting in $L_{2}(B_{r})$ by
\begin{align*}
&T^{(2)}({\delta };\frac{K^{2}}{6})w(p)=\\
& \frac{1}{\pi^2} \int\limits _{B_r} \frac{f(q)d q }
 {(\frac{3 p^2}{4}+1)^{1/4}(p^2+(p,q)+q^2 +1)(\frac{3 q^2}{4}+1)^{1/4}}
  \end{align*}
where $B_{r}=\{p \in T^{3}:|p|<r,\quad r= (\frac{|K|^{2}}{6}
)^{-\frac{1}{2}}\}.$

The equivalence is performed by the unitary dilation

$$U_{r}:\,L_{2}(U_{%
\delta }(0))\rightarrow L_{2}(B_{r}),(U_{r}f)(p)=(\frac{r}{\delta }%
)^{-3/2}f(\frac{\delta }{r}p).$$

Denote by $ \chi_1(\cdot)$ the  characteristic function of $
U_1(0).$  Further, we may replace

 $(\frac{3 p^2}{4}+1)^{-1/4},\,
(\frac{3 q^2}{4}+1)^{-1/4} $ and
$q^2+(p,q)+p^2+1$

by $(\frac{3 p^2}{4})^{-1/4}(1-\chi(p)),\quad
(\frac{3q^2}{4})^{-1/4}(1-\chi(q))$ and $q^2+(p,q)+p^2,$

respectively, since the error
  will be a Hilbert-Schmidt operator continuous up to $K=0.$

Then we get the operator $T^{(2)}(r)$ in $L_2 (U_r(0) \setminus
U_1(0))$ with  the kernel
\[
\frac{2}{\sqrt 3 \pi^2 }\frac{|p|^{-1/2}| q|^{-1/2}}
{q^2+(p,q)+p^2}.
\]
By the dilation $${\bf M}:L_2(U_r(0) \setminus U_1(0))
\longrightarrow L_2((0,{\bf r})\times {\sigma_0}),\,\,\, r=1/2|\log
\frac{|K|^{2}}{6}|, $$
 where
$(M\,f)(x,w)=e^{3x/2}f(e^{ x}w),\, x\in (0,{\bf r}),\, w \in
{\S}^2,$ one sees that the operator $T^{(2)}(r)$ is unitary
equivalent to the integral operator ${\bf S}({\bf r}).$ The
difference of the operators ${\bf S}({\bf r})$ and
$T^{(1)}(\delta,\frac{K^2}{6})$  is compact (up to unitarily
equivalence). Hence Lemma \ref{sobol} yields
$$ \lim\limits_{|K|\to 0}
\frac{n(1,{T}^{(1)}(\delta,\frac{K^2}{6}))} {|log|K||} =U(1).
$$
It is convenient to calculate the coefficient $U(1)$ by means of a
decomposition of the operator
 $\hat{\bf {S}}(y)$ into the orthogonal sum over its invariant
subspaces.

Denote by $L_{l}\subset L_{2}(S^{2})$ the subspace of the  harmonics
of degree $ l=0,1,\cdots. $ It is clear that
$L_{2}(S^{2})=\sum\limits_{l=0}^{\infty }\oplus L_{l},\quad \dim
L_{l}=2l+1.$ Let $P_{l}:\,L_{2}(S^{2})\rightarrow L_{l}$ be the
orthogonal projector onto $L_{l}.$ The kernel of $P_{l}$ is
expressed via the Legendre polynomial $P_{l}(\cdot ):$
\[
  {\ P}_{l}(\xi,\eta)=\frac{2l+1}{4\pi}P_{l}(<\xi,\eta >).
\]
The kernel of $\hat{\bf {S}}(y)$ depends on the scalar product $<\xi
,\eta
>$ only, so that the subspaces $L_{l}$ are invariant for
$\hat{\bf {S}}(y)$ and
\begin{equation*}\label{(5.9)}
  \hat{\bf {S}}(y)=\sum\limits_{l=0}^{\infty}\oplus(\hat{\bf {S}}^{(l)}(y)\otimes{\ P}_{l}),
 \end{equation*}
where $\hat{\mathbf{S}}^{(l)}(y)$ is the multiplication operator by
the number
\begin{equation}\label{(5.10)}
  \hat{S}^{(l)}(y)=2\pi\int_{-1}^{1}P_{l}(t)
  \hat{S}(t;y)dt
\end{equation}
in $L_{l}$ the subspace of the harmonics of degree $l,$ and $P_l(t)$
is a Legendre polynomial. Therefore
$$
  n(\mu,\hat{\bf {S}}(y))=
  \sum\limits_{l=0}^{\infty}(2l+1)n(\mu,\hat{\bf {S}}^{(l)}(y),
  \quad \mu >0.
$$

By \eqref{Stlam} we first calculate \ $\hat{S}^{(0)}(y):$
\begin{equation*}\label{s^0}
\hat{S}^{(0)}(y)=\frac{2}{\sqrt{3}} \int\limits_{-1}^{1}\frac{
\sinh [y(arccos( \frac{1}{2}t)]}{{\sinh (\pi y
)}\sqrt{1-\frac{1}{4}t^{2}}}dt =8\cdot3^{-\frac{1}{2}}\frac{\sinh
\frac{\pi}{6}y} {y\cosh \frac{\pi}{2} y }
\end{equation*}
 It follows from (\ref{sobU}) and  (\ref{(5.10)}) that
\begin{equation*}\label{mes}
{\ U}(1)= \frac{1}{4\pi}\int\limits_{-\infty}^{+\infty} n(1,{\bf
\hat{\bf {S}}}^{(0)}(y))dy=\frac{1}{4\pi}\int\limits_{\hat{S}^{(0)}(\lambda)>1} d\lambda =
\frac{\lambda_0}{2\pi}.
\end{equation*}
In  the same way  proves the second statement of Lemma \ref{main1}.
Lemma \ref{main1} is proven.
\end{proof}

{\bf Proof of Theorem \ref{main}.} Lemmas \ref{comp.pert} ,\ref{b},
\ref{T-cT}, \ref{1Raznost}
 and \ref{main1} yield the proof of Theorem \ref{main}.

 {\bf Proof of Theorem
\ref{asimpZ}  and   \ref{asimpK}.} Let the conditions  of Theorems
\ref{asimpZ}  and   \ref{asimpK} be fulfilled. Then the proof of
Theorems \ref{asimpZ}  and   \ref{asimpK} concerning the
asymptotics
of the number $N(0,z)$ resp.$N(K, 0)$ of eigenvalues follows from
Lemma \ref{b-s} and Theorem \ref{main}.

Now we shall prove the finiteness of $N(K, 0),\,\, K \in
U_\delta^0(0).$

For any $ K \in U_\delta^0(0) $   the kernel
 ${T}^{(1)}(\delta,\frac{K^2}{6};p,q)$ of
  ${T}^{(1)}(\delta,\frac{K^2}{6})$
estimated by
$$|{T}^{(1)}(\delta,\frac{K^2}{6};p,q)|\leq
\frac{C}{|K|^3}, \quad K \in U_\delta^0(0),
 $$
i.e., ${T}^{(1)}(\delta,\frac{K^2}{6})$ belongs to the
Hilbert-Schmidt class.

One concludes from Lemmas \ref{b}, \ref{T-cT} and \ref{1Raznost}
that
 the operator ${\bf T} (K,0)$ splits  in to the  sum of two
 compact(up to
unitarily equivalence) operators
$${\bf T} (K,0)={\bf T}^{(1)}
(K,0)+{\bf \tilde T}(K,0),\,\,\,K \in U_\delta^0(0).$$

Applying Lemma \ref{b-s} and Weyl's inequality
$$n(\lambda_1+\lambda_2,A_1+A_2)\leq
n(\lambda_1,A_1)+n(\lambda_2,A_2)$$  we have
$$
N(K,0)=n(\frac{1}{2},{\bf T}^{(1)} (K,0))+ n(\frac{1}{2},{\bf
\tilde T}(K,0),\,K \in U_\delta^0(0).
$$$\Box$

{\bf Acknowledgement} The authors are  grateful to Prof.Robert A.Minlos
and Prof.Volker Bach for useful discussions.

This work was supported by the DFG 436 USB 113/6 and DFG 436 USB
113/7 projects and the Fundamental Science Foundation of Uzbekistan.
The last two named authors gratefully acknowledge the hospitality of
the Institute of Applied Mathematics and of the IZKS of the
University Bonn.

\end{document}